\documentclass[11pt]{article}
\usepackage{amsmath,amssymb, latexsym,amsfonts}
\usepackage[dvipdf]{graphicx}

\title{Bootstrapping Empirical Processes of Cluster Functionals with Application to Extremograms}
\author{Holger Drees\footnote{University of Hamburg, Department of Mathematics,
 SPST, Bundesstr.\ 55, 20146 Hamburg, Germany;
email: drees@math.uni-hamburg.de}}

\newcommand{\R}{\mathbb{R}}
\newcommand{\N}{\mathbb{N}}
\newcommand{\B}{\mathbb{B}}
\newcommand{\Z}{\mathbb{Z}}

\newcommand{\FF}{\mathcal{F}}

\newcommand{\CC}{\mathcal{C}}
\newcommand{\DD}{\mathcal{D}}

\usepackage[thmmarks,amsmath]{ntheorem}

\theoremstyle{change} \theoremsymbol{$\Box$}
\theoremseparator{\quad}

\theorembodyfont{\em}

\newtheorem{theorem}{Theorem}[section]
\newtheorem{corollary}[theorem]{Corollary}
\newtheorem{remark}[theorem]{Remark}
\newtheorem{proposition}[theorem]{Proposition}

\numberwithin{equation}{section}

\theorembodyfont{\rm}
\newtheorem{example}[theorem]{Example}

\newenvironment{proofof}{\noindent\sc Proof of}{
    \hspace*{\fill} $\Box$ \vspace{2ex} }

\newcommand{\Ind}[1]{{\boldsymbol{1}_{\textstyle\{#1\}}}}

\def\eps{\varepsilon}

\def\Min(#1,#2){#1\wedge #2}
\def\Max(#1,#2){#1\vee #2}
\newcommand{\ceil}[1]{\lceil #1 \rceil}
\newcommand{\floor}[1]{\lfloor #1 \rfloor}

\def\rueck{\noindent\hangafter=1 \hangindent=1.3em}

\begin{document}

\maketitle

\begin{abstract}
   In the extreme value analysis of time series, not only the tail behavior is of interest, but also the serial  dependence plays a crucial role. Drees and Rootz\'{e}n (2010) established limit theorems for a general class of empirical processes of so-called cluster functionals which can be used to analyse various aspects of the extreme value behavior of mixing time series. However, usually the limit distribution is too complex to enable a direct construction of confidence regions. Therefore, we suggest a multiplier block bootstrap analog to the empirical processes of cluster functionals. It is shown that under virtually the same conditions as used by Drees and Rootz\'{e}n (2010), conditionally on the data, the bootstrap processes converge to the same limit distribution. These general results are applied to construct confidence regions for the empirical extremogram introduced by Davis and Mikosch (2009). In a simulation study, the  confidence intervals constructed by our multiplier block bootstrap approach compare favorably to the stationary bootstrap proposed by Davis et al.\ (2012).
\end{abstract}


\section{Introduction}

Time series of observations in environmetrics, (financial) risk
management and other fields often exhibit a non-negligible serial
dependence between extremes. For example, stable areas of low (or
high) pressure may lead to consecutive days of high precipitation
(or high temperature). Likewise, large losses to a financial
investment tend to occur in clusters.

The statistical analysis of the serial dependence structure between
extreme observations is still a challenging task. Yet even if one is
only interested in marginal parameters, like extreme quantiles, it
is crucial to take into account the serial dependence when
assessing the estimation error; see, e.g., Drees (2003) for a
simulation study which demonstrates how misleading confidence
intervals may be if the serial dependence is ignored.

In most applications, no parametric time series model for the
extremal behavior suggests itself. Hence, one should resort to
non-parametric procedures to avoid the risk of an unquantifiable,
but potentially large modeling error. In this context, a general
class of empirical processes that can capture a wide range of
different aspects of the extremal behavior of time series prove a
powerful tool.

To be more concrete, assume that a stationary time series
$(X_t)_{1\le t\le n}$ with values in $E=\R^d$ is observed, from
which we construct $m_n:=\floor{n/r_n}$ blocks
\begin{equation} \label{eq:Ynjdef}
  Y_{n,j} := (X_{n,i})_{(j-1)r_n<i\le jr_n}, \quad 1\le j\le m_n,
\end{equation}
of ``standardized extreme observations'' $X_{n,i}$, $1\le i\le n$. A
typical choice for univariate time series is
\begin{equation} \label{eq:Xnidefuniv}
X_{n,i}:=(X_i-u_n)^+/a_n:= a_n^{-1}(X_i-u_n)1_{\{X_i>u_n\}}
\end{equation}
 for suitable normalizing constants $u_n\in\R$ and $a_n>0$.
Later on, we will use a different notion of extreme
observation in our application to the analysis of the extremogram, for
a multivariate time series.

Denote by $E_\cup := \bigcup_{l\in\N} E^l$ the set of vectors of
arbitrary length with components in $E$, which is equipped with the
$\sigma$-field $\mathbb{E}_\cup$ induced by
    the Borel-$\sigma$-fields on $E^l$, $l\in\N$. Let $\FF$ be a
family of so-called {\em cluster functionals}, i.e.\ functions
$f:(E_\cup,\mathbb{E}_\cup)\to(\R,\B)$ such that $f(0)=0$ and
$f(y_1,\ldots,y_l)=f(0,\ldots,0,y_1,\ldots,y_l,0,\ldots,0)$ for all
$(y_1,\ldots,y_l)\in E_\cup$ where the numbers of coordinates equal
to 0 in the beginning and in the end of the argument on the
right-hand side can be arbitrary. Thus the value of the cluster
functional depends only on the {\em core} of the argument, which is
the smallest subvector of consecutive coordinates that contains all
non-zero values (resp.\ it equals 0 if the argument only consists of
zeros). Then, the pertaining {\em empirical process of cluster
functionals} is defined by
\begin{equation} \label{eq:Zndef}
  Z_n(f) := \frac 1{\sqrt{n v_n}} \sum_{j=1}^{m_n} \big( f(Y_{n,j})-
E f(Y_{n,j})\big), \quad f\in\FF,
\end{equation}
with $v_n:= P\{X_{n,1}\ne 0\}$. Drees and Rootz\'{e}n (2010) established
sufficient conditions for $Z_n$ to converge to a Gaussian process in
the space $\ell^\infty(\FF)$ of bounded functions on $\FF$. The
following theorem summarizes their main results; the conditions are
recalled in the appendix.
\begin{theorem} \label{theo:DRmain}
  \begin{enumerate}
    \item If the conditions (B1), (B2) and (C1)--(C3) are fulfilled,
    the finite-dimen\-sio\-nal marginal distributions (fidis) of the
    empirical process $Z_n$ converge to the pertaining fidis of a
    Gaussian process $Z$ with covariance function $c$ (defined in
    (C3)).
    \item Under the conditions (B1), (B2) and (D1)--(D4) the
    empirical process $Z_n$ is asymptotically tight in
    $\ell^\infty(\FF)$. If, in addition, the conditions (C1)--(C3)
    are met, then $Z_n$ weakly converges to $Z$.
    \item If the assumptions (B1), (B2), (D1), (D2'), (D3) and (D5)
    are satisfied and, in addition, (D6) (or the more restrictive
    condition (D6')) holds, then $Z_n$ is asymptotically
    equicontinuous. Hence, $Z_n$ weakly converges to $Z$ in
    $\ell^\infty(\FF)$ if also the conditions (C1)--(C3) hold.
  \end{enumerate}
\end{theorem}

 For certain types of families $\FF$ of cluster functionals, Drees and
Rootz\'{e}n (2010) also gave sets of conditions that are sufficient for
$(Z_n(f))_{f\in\FF}$ to converge and easier to verify than the
abstract conditions listed in the appendix.

We will demonstrate their usefulness by improving on limit results
on an empirical version of the so-called extremogram introduced by
Davis and Mikosch (2009) in the framework of time series with
regularly varying marginals. To be more precise, assume that
$(X_t)_{t\in\Z}$ is a stationary $\R^d$-valued time series such that
for all $h\in\N$ the vector $(X_0,X_h)\in\R^{2d}$ is regularly
varying. Recall that a random vector $W\in\R^l$ is regularly varying
if there exists a non-null  measure $\nu$ on
$\R^l\setminus\{0\}$
such that
$$ \frac{P\{W\in xB\}}{P\{\|W\|>x\}}\;\longrightarrow\;
\nu(B)<\infty
$$
for all $\nu$-continuity sets $B\in\B^l$ that are bounded away from
the origin 0. Note that, while this definition of regular variation
does not depend on the choice of the norm $\|\cdot\|$, the specific
form of the limiting measure $\nu$ does. In  any case, the limiting measure is homogeneous of order $-\alpha$ for some $\alpha>0$, the so-called index of regular variation.

Then, with $F_{\|X\|}^\leftarrow$ denoting the quantile function of
$\|X_0\|$ and  $a_n:= F_{\|X\|}^\leftarrow(1-1/n)\to\infty$, to each
lag $h\in\N$ there exists a measure $\nu_{(0,h)}$ on
$\R^{2d}\setminus\{0\}$ such that
\begin{equation}  \label{eq:nuhdef}
  nP\{a_n^{-1}(X_0,X_h)\in B\}\;\longrightarrow\;
\nu_{(0,h)}(B)
\end{equation}
for all $\nu_{(0,h)}$-continuity sets $B\in\B^{2d}$ bounded away from the origin.
In particular, for all $A,B\in\B^d$ bounded away from 0 such that
$\nu_h(\partial(A\times B))=0=\nu_h(\partial(A\times \R^d))$ and
$\nu_h(A\times \R^d)>0$ one has
$$ P(X_h\in a_nB\mid X_0\in a_n A) = \frac{P\{a_n^{-1}(X_0,X_h)\in
A\times B\}}{P\{a_n^{-1}X_0\in A\}} \;\longrightarrow\;
\frac{\nu_{(0,h)}(A\times B)}{\nu_{(0,h)}(A\times R^d)} =:
\rho_{A,B}(h).
$$
Davis and Mikosch (2009) called $\rho_{A,B}$ (as a function of $h$)
the {\em extremogram} of $(X_t)_{t\in\Z}$ (pertaining to $A,B$). It
is worth mentioning that the extremogram is closely related to the
concept of tail processes introduced by Basrak and Segers (2008).

Based on the observations $X_1,\ldots, X_n$, they proposed the
following empirical counterpart as an estimator of $\rho_{A,B}(h)$:
\begin{equation}  \label{eq:empextremodef1}
  \hat \rho_{A,B}(h) := \frac{\sum_{i=1}^{n-h} \Ind{X_i\in a_k A,
  X_{i+h}\in a_kB}}{\sum_{i=1}^n \Ind{X_i\in a_k A}}.
\end{equation}
Here $k=k_n$ is a sequence that tends to $\infty$ at a slower rate
than $n$ so that $a_k\to \infty$ at a slower rate than $a_n$, and
thus the number of extreme observations used for estimation tends to
$\infty$. Under suitable conditions, $(\hat
\rho_{A,B}(h))_{h\in\{0,\ldots,h_0\}}$ is asymptotically normal (see
Davis and Mikosch, 2009, Corollary 3.4).

This result has two serious drawbacks. First, usually, the
normalizing constants $a_k$ are unknown and must hence be replaced
with an empirical counterpart, like, e.g., the $\floor{n/k}+1$
largest observed norm:
\begin{equation} \label{eq:amhatdef}
  \hat a_k := \hat a_{k,n} := \|X\|_{n-\floor{n/k}:n}.
\end{equation}
It is not obvious whether this modification influences the
asymptotic behavior of the empirical extremogram.

Secondly, the extremogram for a fixed pair of sets $A$ and $B$
conveys limited information on the extremal dependence structure, in
particular in a multivariate setting, i.e.\ if $d>1$. To get a
fuller picture, one should consider the extremogram for a whole
family of sets simultaneously. For example, in the case $d=1$, Drees et al.\ (2015) considered rays $(-\infty,-x)$ and $(x,\infty)$ for all $x>0$ simultaneously. However, the techniques used by Davis
and Mikosch (2009) are not applicable to infinite families of sets.

We will show that both problems can be neatly solved using the
theory of empirical processes of cluster functionals. Indeed, if the
families of sets $A$ and $B$ are suitably chosen and the bias of
$\hat \rho_{A,B}(h) $ is asymptotically negligible, then the
asymptotic normality of the empirical extremogram with estimated
normalizing sequence $\hat a_k$ follows immediately.

If one wants to construct confidence regions using this limit
theorem, then estimators of the limiting covariance structure are
needed. Since the direct estimation does not look promising, Davis
et al.\ (2012) proposed to use a so-called stationary bootstrap
instead. Here we follow a somewhat different approach. First, in the
general setting considered by Drees and Rootz\'{e}n (2010), it is shown
that the convergence of a multiplier block bootstrap version of the
empirical process of cluster functional conditionally given the data
follows under the same conditions as the convergence of $Z_n$
itself. From this powerful result it is easily concluded that a
multiplier block bootstrap version can be used to construct
confidence regions for the extremogram.

Though in the present paper we focus on the extremogram as one
possible measure for the extremal dependence structure of the time
series, the same approach using empirical processes of cluster
functionals can be used in a much wider context. For example, Drees
(2011) analyzed block estimators of the so-called extremal index of
absolutely regular time series using empirical processes of cluster
functionals and suggested a bias corrected version thereof.

The paper is organized as follows. In Section 2 we introduce
multiplier block bootstrap versions of the empirical process $Z_n$.
Moreover, we give sufficient conditions under which, in probability
conditional on the data, this bootstrap processes weakly converge to
the same limiting process as $Z_n$. In Section 3, it is demonstrated
that the theory developed by Drees and Rootz\'{e}n (2010) yields limit theorems for the empirical extremogram with
estimated normalizing sequence uniformly over suitable families of
sets. In the same setup, a bootstrap result easily follows from the
general theory developed in Section 2. The results of a small simulation study are reported in Section 3. All proofs are postponed to
Section \ref{proofs}.

Throughout the paper, we will use the notation $x^{(k)}$ for the
vector $(x_1, \ldots,x_k)$  made up by the first $k$ components in
the vector $x$, if $x$ has at least $k$ components, and otherwise
$x^{(k)}=x$. The maximum norm of a vector $x\in\R^l$ for some
$l\in\N$ is denoted by $\|x\|$. We omit indices of random variables
to denote a generic random variable with the same distribution; for
example, $\xi$ is a generic random variable with the same
distribution as $\xi_j$ and $Y_n$ is a generic random vector with
the same distribution as $Y_{n,j}$.

\section{Multiplier processes}  \label{sect:bootstrap}

In what follows, $(X_{n,i})_{1\le i\le n, n\in\N}$ is a row-wise
stationary triangular scheme of $E=\R^d$-valued random vectors.
Usually these vectors are derived from some fixed stationary time
series $(X_t)_{t\in\Z}$ by a transformation which depends on the
stage $n$ and which sets all but the ``extreme'' observations to 0
in such a way that the probability that a transformed observation is
non-zero tends to 0 as $n\to\infty$. For univariate time series,
often definition \eqref{eq:Xnidefuniv} is used. In our application
to the empirical extremogram instead we define
\begin{equation}  \label{eq:Xnidefextremo}
  X_{n,i}^{(h,\tilde h)} := a_k^{-1} \big(X_i\Ind{X_i\not\in (-\infty,a_kx_*)^d},X_{i+h}\Ind{X_{i+h}\not\in (-\infty,a_kx_*)^d},X_{i+\tilde h}\Ind{X_{i+\tilde h}\not\in (-\infty,a_kx_*)^d}\big)
\end{equation}
for  some $x_*>0$ and $h,\tilde h\in\N_0$.

According to Theorem \ref{theo:DRmain}, under suitable conditions,
the empirical process $Z_n$ of cluster functionals converge to a
Gaussian process $Z$ with covariance function $c$, which is defined
in (C3) as the limit of the covariance function of the cluster
functionals applied to a block $Y_n$ of $r_n$ consecutive
``standardized extremes'' $X_{n,i}$. One may try to estimate this
covariance function by an empirical covariance, but since most of
the blocks $Y_{n,j}$ defined in \eqref{eq:Ynjdef} equal 0, a
bootstrap approach seems more promising.

Because the processes are defined via functionals applied to whole
blocks $Y_{n,j}$ of ``standardized extremes'', it suggests itself to
use some block bootstrap. More precisely, we consider the following
two versions of multiplier block bootstrap processes:
\begin{eqnarray}  \label{eq:multiplierdef}
  Z_{n,\xi}(f) & := & \frac 1{\sqrt{n v_n}} \sum_{j=1}^{m_n} \xi_j \big( f(Y_{n,j})-
E f(Y_{n,j})\big),\\
 Z_{n,\xi}^*(f) & := & \frac 1{\sqrt{n v_n}} \sum_{j=1}^{m_n} \xi_j \big( f(Y_{n,j})-
\overline{f(Y_n)}\big), \quad f\in\FF, \label{eq:bootstrapdef}
\end{eqnarray}
where $\overline{f(Y_n)} := m_n^{-1} \sum_{j=1}^{m_n} f(Y_{n,j})$
and $\xi_j, j\in\N$, are i.i.d.\ random variables with $E(\xi_j)=0$
and $Var(\xi_j)=1$ independent of $(X_{n,i})_{1\le i\le n, n\in\N}$.
Note that in the definition of the multiplier process $Z_{n,\xi}$
expectations $E f(Y_n)$ are used which are usually unknown to the
statistician. Hence, in some applications, it may be useful to
replace them with the estimators $\overline{f(Y_n)}$, which leads to
the bootstrap processes $Z_{n,\xi}^*$.

Our main goal is to prove weak convergence of $Z_{n,\xi}$ and
$Z_{n,\xi}^*$ to $Z$ in probability, conditionally on the data. To
this end, as usual, we  metrize weak convergence in
$\ell^\infty(\FF)$ using the bounded Lipschitz metric on the space
of probability measures on $\ell^\infty(\FF)$. That is, for two
probability measures $Q_1$ and $Q_2$ we define
$$ d_{BL(\ell^\infty(\FF))}(Q_1,Q_2) := \sup_{g\in
BL_1(\ell^\infty(\FF))} \Big| \int g\, dQ_1-\int g\, dQ_2\Big|,
$$
where
\begin{eqnarray*}
 BL_1(\ell^\infty(\FF)) & := & \big\{
g:\ell^\infty(\FF)\to \R\mid \|g\|_\infty := \sup_{z\in
\ell^\infty(\FF)}|g(z)|\le 1,\\
& & \hspace*{0.5cm} |g(z_1)-g(z_2)|\le \|z_1-z_2\|_\FF :=
\sup_{f\in\FF}|z_1(f)-z_2(f)| \text{ for all }
z_1,z_2\in\ell^\infty(\FF)\big\}.
\end{eqnarray*}
 Likewise, for the convergence of the fidis, we use the distance
$$ d_{BL(\R^l)}(Q_1,Q_2) := \sup_{g\in
BL_1(\R^l)} \Big| \int g\, dQ_1-\int g\, dQ_2\Big|,
$$
between two probability measures $Q_1$ and $Q_2$ on $\R^l$, where
$$ BL_1(\R^l) := \big\{ g:\R^l\to \R\mid
\sup_{v\in \R^l}|g(v)|\le 1, |g(v_1)-g(v_2)|\le \|v_1-v_2\| \text{
for all } v_1,v_2\in\R^l\big\}.
$$

By $E_\xi$ (resp.\ $E_\xi^*$) we denote the (outer) expectation with
respect to $(\xi_j)_{j\in\N}$, i.e.\linebreak
$E_\xi\big(f(\xi_1,\ldots,\xi_{m_n},X_{n,1},\ldots,
X_{n,n})\big)=E\big(f(\xi_1,\ldots,\xi_{m_n},X_{n,1},\ldots,
X_{n,n}) \mid X_{n,1},\ldots, X_{n,n}\big)$ is the expectation of
the function conditionally on the observations. Likewise, we denote by $P_\xi$ the
probability measure w.r.t.\ $(\xi_j)_{j\in\N}$. (Cf.\ Kosorok, 2003,
for a precise definition using a special construction of probability
spaces.)

Our first result shows that the asymptotic behavior of the fidis of
$Z_{n,\xi}$, conditionally on the data, is the same as the
(unconditional) behavior of the fidis of $Z_n$.
\begin{theorem}  \label{theo:bootfidis}
  Under the conditions (B1), (B2) and (C1)--(C3) one has for all
  $f_1,\ldots, f_l\in\FF$
  \begin{equation}   \label{eq:bootfidis}
   \sup_{g\in BL_1(\R^l)} \Big| E_\xi g\big((Z_{n,\xi}(f_k))_{1\le
  k\le l}\big) -E g\big((Z(f_k))_{1\le
  k\le l}\big) \Big| \;\longrightarrow\; 0
  \end{equation}
  in probability.
\end{theorem}
Since the supremum in \eqref{eq:bootfidis} is bounded by 2, it readily follows that
\begin{eqnarray*}
  \lefteqn{\sup_{g\in BL_1(\R^l)} \Big| E g\big((Z_{n,\xi}(f_k))_{1\le
  k\le l}\big) -E g\big((Z(f_k))_{1\le
  k\le l}\big) \Big|}\\
  & \le & E \sup_{g\in BL_1(\R^l)} \Big| E_\xi g\big((Z_{n,\xi}(f_k))_{1\le
  k\le l}\big) -E g\big((Z(f_k))_{1\le
  k\le l}\big) \Big| \;\longrightarrow\; 0,
  \end{eqnarray*}
that is, the (unconditional) weak convergence of the fidis of
$Z_{n,\xi}=(Z_{n,\xi}(f))_{f\in\FF}$ to the corresponding fidis of
$Z$.

Following the ideas developed by Kosorok (2003), the following
result establishes the asymptotic tightness of $Z_{n,\xi}$ under a
bracketing entropy condition, and thus also the weak convergence of
$Z_{n,\xi}$ under the same conditions as the convergence of the
original empirical process in Theorem \ref{theo:DRmain}(ii).
\begin{proposition}  \label{prop:boottight}
  Suppose that the conditions (B1), (B2), (D1), (D3) and (D4) hold
  and
  \begin{enumerate}
    \item (D2) holds and $\xi$ is bounded, or
    \item (D2') holds and $E^*(F^2(Y_n))=O(r_n v_n)$.
  \end{enumerate}
  Then $Z_{n,\xi}$ is asymptotically tight in $l^\infty(\FF)$. Hence
  it converges to $Z$ if, in addition, the conditions (C1)--(C3) are
  met.
\end{proposition}

Now a modification of the arguments given in the proof of Theorem 2
of Kosorok (2003) yields the desired convergence result for the
multiplier process conditionally on the data.
\begin{theorem}  \label{theo:bootconv1}
  If condition (D3) and convergence \eqref{eq:bootfidis} hold and
  $Z_{n,\xi}$ weakly converges to $Z$, then
  \begin{equation}  \label{eq:bootconv1}
    \sup_{g\in BL_1(\ell^\infty(\FF))} \big| E_\xi g(Z_{n,\xi})-
    Eg(Z)\big| \;\longrightarrow\; 0
  \end{equation}
  in outer probability.
\end{theorem}
A combination of this result with Theorem \ref{theo:bootfidis} and
Proposition \ref{prop:boottight} leads to
\begin{corollary} \label{corol:bootconv1}
  If the conditions (B1), (B2), (C1)-(C3) and (D1)-(D4) are
  satisfied and $\xi$ is bounded, then convergence \eqref{eq:bootconv1} holds.
\end{corollary}

According to Theorem \ref{theo:bootconv1}, under (D3) the weak
convergence of the multiplier process $Z_{n,\xi}$ to $Z$
conditionally on the data follows from the weak convergence of the
fidis conditionally on the data and the (unconditional) convergence
of $Z_{n,\xi}$ to $Z$. The latter assertion may also be derived by
establishing the asymptotic equicontinuity of $Z_{n,\xi}$ using a
metric entropy condition (instead of verifying tightness using a
bracketing entropy condition as in Proposition
\ref{prop:boottight}).
\begin{proposition}  \label{prop:bootequi}
  Suppose that the conditions (B1), (B2), (D1), (D2'), (D3)
  and\\[1ex]
  {\bf (D5')}\;\;\;  \parbox[t]{14.9cm}{For all $\delta>0, n\in\N, (e_i)_{1\le i\le
   \floor{m_n/2}} \in \{-1,0,1\}^{\floor{m_n/2}}$ and $k\in\{1,2\}$ the
   map $\sup_{f,g\in\FF, \rho(f,g)<\delta}$
   $\sum_{j=1}^{\floor{m_n/2}} e_j\big(\xi_j(
   f(Y_{n,j}^*)-g(Y_{n,j}^*))\big)^k$
      is measurable}\\[2ex]
  are fulfilled and
  \begin{enumerate}
    \item (D6) holds and $\xi$ is bounded, or
    \item (D6') holds.
  \end{enumerate}
  Then $Z_{n,\xi}$ is asymptotically equicontinuous. Hence,
  it converges to $Z$ if, in addition, the conditions (C1)--(C3) are
  met.
\end{proposition}

Using Theorem \ref{theo:bootconv1} and Corollary 2.6.12 of van der
Vaart and Wellner (1996), we obtain as an immediate consequence
\begin{corollary}  \label{corol:bootconv2}
  If the conditions (B1), (B2), (C1)-(C3), (D1), (D2'), (D3) and
  (D5') are met, if $F$ is measurable with $E(F^2(Y_n))=O(r_nv_n)$ and $\FF$ is a VC-hull
  class, then convergence \eqref{eq:bootconv1} holds.
\end{corollary}

To sum up, we have shown that, roughly under the same conditions as
used in Theorem \ref{theo:DRmain}, the multiplier process
$Z_{n,\xi}$ shows the same asymptotic behavior conditionally on the
data as the empirical process $Z_n$ unconditionally. The following
result gives conditions under which the convergence of $Z_{n,\xi}$
implies the convergence of the bootstrap process $Z_{n,\xi}^*$
conditionally on the data.
\begin{corollary}  \label{corol:Znxistarconv}
  If convergence \eqref{eq:bootfidis} of the fidis of $Z_{n,\xi}$
   holds conditionally on the data, condition (D3) is satisfied and
   $Z_n\to Z$ and $Z_{n,\xi}\to Z$ weakly, then
   \begin{equation}  \label{eq:ZnxiZnxistarapprox}
    E_\xi \sup_{f\in\FF}
   |Z_{n,\xi}^*(f)-Z_{n,\xi}(f)|\;\longrightarrow\; 0
   \end{equation}
   in outer probability, $Z_{n,\xi}^*\to Z$ weakly and
   \begin{equation}  \label{eq:bootconv2}
    \sup_{g\in BL_1(\ell^\infty(\FF))} \big| E_\xi g(Z_{n,\xi}^*)-
    Eg(Z)\big| \;\longrightarrow\; 0
  \end{equation}
  in outer probability. In particular, these assertions hold under
  the conditions of Corollary \ref{corol:bootconv1} and under the
  assumptions  of Corollary \ref{corol:bootconv2}.
\end{corollary}

\begin{remark} \label{rem:estnorm}
  Note that also the normalizing factor $(nv_n)^{-1/2}$ in the
  definition of $Z_{n,\xi}^*$ may be unknown. In most
  applications of multiplier processes, though, this is not
  problematic, because this factor is not needed to construct
  confidence regions. Nevertheless, it is noteworthy that assertion
  \eqref{eq:bootconv2} remains valid if $v_n$ is replaced with some
  estimator $\hat v_n$ that is consistent in the sense that  $\hat
  v_n/v_n\to 1$ in probability.
\end{remark}

For specific types of cluster functionals, Drees and Rootz\'{e}n (2010)
gave simpler sufficient conditions for the convergence of the
corresponding empirical process which carry over to the multiplier
processes considered here. In the next section we will use the
conditions of Corollary 3.6 of that paper, which deals with
so-called generalized tail array sums, i.e.\ empirical processes
with functionals of the form $f_\phi(y_1,\ldots,y_r)=\sum_{i=1}^r
\phi(y_i)$ for functions $\phi:(E,\mathbb{E})\to(\R,\B)$ such that
$\phi(0)=0$.

\section{Processes of Extremograms}

In this section we employ the general theory to analyze the
asymptotic behavior of the empirical extremogram $\hat\rho_{n,A,B}$,
a version with empirical normalization and a bootstrap version
thereof, uniformly over suitable families of sets $A$ and $B$ and
over lags $h\in\{0,\ldots,h_0\}$ for some fixed $h_0\in\N$.
Throughout this section we are only interested in the behavior for
vectors with at least one large component. We thus consider families
$\CC$  of pairs of measurable subsets of $\R^d$ such that
$$ x_* := \inf_{(A,B)\in\CC} \inf_{x\in A} \max_{1\le j\le d}
x_j>0,
$$
i.e.\ $A\subset \R^d\setminus (-\infty,x_*)^d$ for all
$(A,B)\in\CC$. However, the results below can be generalized to
families of sets that are uniformly bounded away from 0 so that
$\inf_{(A,B)\in\CC}\inf_{x\in A} \max_{1\le j\le d} |x_j|>0 $. For
the sake of notational simplicity, we assume that $n+h_0$ (instead
of $n$) $\R^d$-valued random vectors $X_1,\ldots,X_{n+h_0}$ are
observed.
\begin{remark} \label{rem:cones}
  To keep the presentation simple, we will assume that $X_0$ is
  regularly varying on the full cone
  $\R^d\setminus\{0\}$ with a limiting measure $\nu_0$
  which is not concentrated on $(-\infty,0]^d$; see Theorem
  \ref{theo:extremo1} below. This assumption could be weakened to
  the regular variation on the cone $\R^d\setminus (-\infty,0]^d$
  defined in the spirit of Das et al.\ (2013), i.e.\ there exists a
  normalizing sequence $\tilde a_n>0$ and a measure $\tilde \nu_0$ such that
  $$ nP\{ X_0/\tilde a_n\in B\} \;\longrightarrow\; \tilde \nu_0(B) $$
  for all $\tilde\nu_0$-continuity  sets $B\in\B$ bounded away from
  $(-\infty,0]^d$, where the limit has to be finite. Here one may
  choose $\tilde a_n$ as the $(1-1/n)$-quantile of $\max_{1\le j\le
  d} X_{0,j}$. Under the slightly more restrictive assumption used
  in the results below, one has
  $$ \frac{P\{\max_{1\le j\le d} X_{0,j}>u\}}{P\{\|X_0\|>u\}}
  \;\longrightarrow\;\nu_0\big(\R^d\setminus(-\infty,1]^d\big)
  $$
  as $u\to\infty$, and hence $\tilde a_n\sim a_n
  \big(\nu_0\big(\R^d\setminus(-\infty,1]^d\big)\big)^{1/\alpha}$ and
  $\tilde\nu_0=\nu_0/\nu_0\big(\R^d\setminus(-\infty,1]^d\big)$,
  where $-\alpha$ is the degree of homogeneity of $\nu_0$, i.e.\
  $\nu_0(\lambda B)=\lambda^{-\alpha}\nu_0(B)$.
\end{remark}

For some intermediate sequence $k=k_n$ (i.e.\ $k_n\to\infty$,
$k_n/n\to 0$), we define the empirical extremogram to the sets $A$
and $B$ and lag $h$ as
$$ \hat \rho_{n,A,B}(h) := \frac{\textstyle \sum_{i=1}^n
1_{A\times B}(X_i/a_k,X_{i+h}/a_k)}{\textstyle \sum_{i=1}^n
1_{A}(X_i/a_k)}.
$$
Note that this is a slight modification of the definition given by
Davis and Mikosch (2009) in that we do not use the maximal number of
summands in the denominator. However, it is easily seen that all
results given below carry over to the original definition.

The uniform asymptotic behavior of the empirical extremogram will
easily follow from that of the stochastic process
\begin{eqnarray*}
 \tilde Z_n(h,A,B) & := & \frac 1{\sqrt{nv_n}} \sum_{i=1}^n \Big( 1_{A\times
B}(X_i/a_k,X_{i+h}/a_k)-P\{X_i\in a_kA,X_{i+h}\in a_kB\}\Big),
\end{eqnarray*}
$h\in\{0,\ldots,h_0\}, (A,B)\in\CC,$ with
$$ v_n:=  P\{X_0\not\in (-\infty,a_kx_*)^d \}.
$$
This process, in turn, can be analyzed using the theory for
empirical processes of cluster functionals developed by Drees and
Rootz\'{e}n (2010). In order to use conditions on the joint distribution
of the $X_t$ as weak as possible, it is useful to consider such
processes indexed by $(A,B)\in\CC$ and just two lags $h,\tilde
h\in\{0,\ldots,h_0\}$. Let
\begin{eqnarray*}
  \tilde X_{n,i} & := & \frac{X_i}{a_k} 1_{\R^d\setminus
  (-\infty,x_*)^d}\Big(\frac{X_i}{a_k}\Big), \quad 1\le i\le
  n+h_0,\\
  X_{n,i}^{(h,\tilde h)} & := & (\tilde X_{n,i},\tilde X_{n,i+h},\tilde X_{n,i+\tilde h}), \quad 1\le i\le
  n,\\
 Y_{n,j}^{(h,\tilde h)} & := & (X_{n,i}^{(h,\tilde
 h)})_{(j-1)r_n<i\le jr_n}, \quad 1\le j\le m_n,\\
   v_n^{(h,\tilde h)} & := & P\{X_{n,i}^{(h,\tilde h)}\ne 0\} =
   P\{(X_0,X_h,X_{\tilde h})\not\in (-\infty,x_*)^{3d}\},\\
  \DD & := & \{A\times B\times \R^d, A\times\R^d\times B\mid
  (A,B)\in \CC\},\\
  f_D(y_1,\ldots, y_r) & := & \sum_{i=1}^r 1_D(y_i), \quad
  y_i\in\R^{3d},\quad  D\in\DD,\\
  \FF & := & \{f_D\mid D\in\DD\}, \quad \text{and}\\
  Z_n^{(h,\tilde h)}(f_D) & := & \frac 1{\sqrt{n v_n^{(h,\tilde h)}}}
  \sum_{j=1}^{m_n} \big( f_D(Y_{n,j}^{(h,\tilde h)})- E f_D(Y_{n,j}^{(h,\tilde
  h)})\big)\\
   & = & \frac 1{\sqrt{n v_n^{(h,\tilde h)}}}
  \sum_{i=1}^{m_nr_n} \big( 1_D(X_{n,i}^{(h,\tilde h)})- P\{X_{n,i}^{(h,\tilde
  h)}\in D\}\big), \quad D\in\DD.
 \end{eqnarray*}

Note that, for $n=m_nr_n$, we have $\tilde
Z_n(h,A,B)=(v_n^{(h,\tilde h)}/v_n)^{1/2} Z_n^{(h,\tilde
h)}(f_{A\times B\times \R^d})$ and $\tilde Z_n(\tilde
h,A,B)=(v_n^{(h,\tilde h)}/v_n)^{1/2} Z_n^{(h,\tilde h)}(f_{A\times
\R^d\times B})$; under the conditions of Theorem \ref{theo:extremo1}
the difference between these processes is asymptotically negligible
even if $m_nr_n<n$.

Using Corollary 3.6 of Drees and Rootz\'{e}n (2010) and Drees and Rootz\'{e}n (2015), we
obtain the following set of sufficient conditions for the
convergence of $\tilde Z_n$.
\begin{theorem} \label{theo:extremo1}
  Suppose that all four-dimensional marginal distributions of the
  stationary time series $(X_t)_{t\in\N_0}$ are regularly varying, i.e.\ for all index vectors $I\in\N_0^l$ of dimension $l\le 4$ there exists a measure $\nu_I$ such that
  \begin{equation} \label{eq:fourdimregvar}
   n P\{a_n^{-1}X_I\in B\} \;\longrightarrow\,
  \nu_I(B)<\infty
  \end{equation}
  for all Borel sets $B$ bounded away from $0\in\R^{ld}$, and that $\nu_0(\R^d\setminus(-\infty,x^*)^d)>0$.
   In addition, assume  that
  the conditions (B1), (B2) and ($\widetilde{\text{B3}}$) are fulfilled, and
    $r_n=o(\sqrt{nv_n})$. Finally, assume that there exists a bounded semi-metric $\bar\varrho$
   on $\CC$ such that $\CC$ is totally bounded w.r.t.\
  $\bar\varrho$, and a function $u:(0,\infty)\to (0,\infty)$ such that
  $\lim_{t\downarrow 0} u(t)=0$ and
  \begin{equation}   \label{eq:clustermoment}
   E\Big(\sum_{i=1}^{r_n} 1_{(A\times B)\Delta (\tilde A\times\tilde
  B)}(X_i/a_k,X_{i+h}/a_k)\Big)^2 \le u\big(
  \bar\varrho\big((A,B),(\tilde A,\tilde B)\big)\big) r_n v_n
  \end{equation}
  for all $(A,B),(\tilde A,\tilde B)\in\CC$, $h\in\{0,\ldots,h_0\}$,
   and that the conditions (D5) and (D6) hold for
   $\varrho\big(f_D,f_{\tilde D}\big):= \bar\varrho\big((A,B),(\tilde A,\tilde
   B)\big)$ if $D=A\times B\times \R^d$, $\tilde D=\tilde A\times \tilde B\times
   \R^d$, or $D=A\times \R^d\times B$, $\tilde D=\tilde A\times \R^d\times \tilde B$,
   and $\varrho\big(f_D,f_{\tilde D}\big):= L$ else for some
   sufficiently large constant $L>1$.
   (Here $C_1\Delta C_2$ denotes the
symmetric difference of the two sets $C_1$ and $C_2$.)

    Then
    $\tilde Z_n$ converges weakly to
  a  Gaussian process $\tilde Z$ with covariance function
   $$\tilde  c\big((h,A,B),(\tilde h,\tilde A,\tilde B)\big):=
   \sum_{i=-\infty}^\infty \frac{\nu_{(0,h,i,i+\tilde h)}(A\times B\times \tilde A\times \tilde B)}{
   \nu_{0}\big(\R^d\setminus (-\infty,x_*)^d\big)}<\infty.
   $$
   { }
\end{theorem}
Observe that in \eqref{eq:fourdimregvar} necessarily the following
consistency condition holds: for vectors $I_0=(i_j)_{1\le j\le l}$
and $I=(i_j)_{1\le j\le 4}$ of indices and $\nu_{I_0}$-continuity
sets $A\in \B^{ld}$ bounded away from the origin one has
$\nu_{I_0}(A)=\nu_I(A\times\R^{(4-l)d})$.

Usually the moment condition \eqref{eq:clustermoment} and  the
entropy condition (D6) are most difficult to verify. The proof
of Theorem \ref{theo:extremo1}  shows that the process $\tilde Z_n$
indexed by $\tilde \FF:=\{(h,A,B)\mid h\in\{0,\ldots,h_0\},
(A,B)\in\CC\}$ is asymptotically tight if and only if the empirical
processes $Z_n^{(h,\tilde h)}$ indexed by $\{f_{A\times B\times
\R^d}\mid (A,B)\in\CC\}$ resp.\ $\{f_{A\times \R^d\times B}\mid
(A,B)\in\CC\}$ are asymptotically tight for all $h,\tilde
h\in\{0,\ldots,h_0\}$. Thus we may replace condition (D6) by the
assumption that these families are VC-subgraph class of functions,
which in turn is equivalent to the assumption that
\begin{equation} \label{eq:defFbar}
\bar\FF:=\{\bar
f_{A\times B}\mid (A,B)\in\CC\} \quad \text{with} \quad \bar f_D(y_1,\ldots, y_r)  :=
\sum_{i=1}^r 1_D(y_i) \;\; \text{for} \;\; y_i\in\R^{2d}, 1\le i\le r,
\end{equation}
 is a VC-subgraph class of
functions. Likewise, one may divide the family $\CC$ into a finite
number of subfamilies $\CC_j$ and check that $\bar\FF_j:=\{\bar
f_{A\times B}\mid (A,B)\in\CC_j\}$ is a VC-subgraph class of
functions.

For applications to the asymptotic analysis of empirical
extremograms, we shall consider families $\CC$ such that for
$(A,B)\in\CC$ also $(A,\R^d)$ belongs to $\CC$. The following simple
example exhibits another closedness property of $\CC$ which is
important to prove convergence of the empirical extremogram with
estimated normalizing constant.
\begin{example} \label{ex:multiplesets}\rm
  Fix some $\lambda_0>0$ and measurable sets $A_0,B_0\subset \R^d$
  bounded away from 0 such that $x\in A_0$ implies $\lambda x\in A_0$
  for all $\lambda>1$ and likewise for $x\in B_0$. (In particular, one
  may choose a set $A_0\subset [0,\infty)^d\setminus[0,1]^d$ such
  that $x\in A_0$ and $y\ge x$ imply $y\in A_0$.) Then, for
  $\CC_1:=\{\lambda(A_0,B_0)\mid \lambda> \lambda_0\}$, the family
  $\bar\FF_1$ is a VC-subgraph class of functions. To see this, note
  that $\bar f_{\lambda (A\times  B)} \le f_{\tilde \lambda
  (A\times B)}$ if $\lambda> \tilde\lambda$, i.e.\ the functions are linearly ordered. Hence no set of size 2 can be
  shattered by the subgraphs of $\bar\FF_1$. Likewise, the family
  $\bar\FF_2$ pertaining to $\CC_2:=\{\lambda(A_0,\R^d)\mid \lambda>
  \lambda_0\}$  is a VC-subgraph class of functions.

  Condition \eqref{eq:clustermoment} can be reformulated as follows.
  There exists a semi-metric $\tilde\varrho$ on $[\lambda_0,\infty)$
  such that $[\lambda_0,\infty)$ is totally bounded w.r.t.\
  $\tilde\varrho$ and $E\big(\sum_{i=1}^{r_n} 1_{(\lambda(A_0\times
  B_0))\setminus (\tilde\lambda(A_0\times
  B_0))}(X_i/a_k,X_{i+h}/a_k)\big)^2\le
  u(\tilde\varrho(\lambda,\tilde\lambda))r_nv_n$ and  $E\big(\sum_{i=1}^{r_n}
  1_{(\lambda A_0)\setminus (\tilde\lambda A_0)}(X_i/a_k)\big)^2\le
  u(\tilde\varrho(\lambda,\tilde\lambda))r_nv_n$ hold for all
  $\lambda_0< \lambda<\tilde\lambda$ and all $n\in\N$.
\end{example}

The families of sets $A$ and $B$  most widely discussed in the
literature are sets of upper right orthants $(x,\infty)$ and
complements $\R^d\setminus (-\infty,x]$ of lower left orthants.
\begin{example} \label{ex:orthants} \rm
  Consider the family $\CC_1:=\big\{ \big((x_A,\infty),(x_B,\infty)\big) \mid
  x_A,x_B \not\in (-\infty,x_*]^d\big\}$ of pairs of upper right
  orthants bounded away from the origin. Then condition (D6) holds for
   $\CC:=\CC_1\cup\big\{\big((x_A,\infty),\R^d\big)|x_A\not\in (-\infty,x_*]^d\big\}$ if condition (B1) is satisfied and
   \begin{equation} \label{eq:extremomomcond}
     E\Big( \sum_{i=1}^{r_n} 1_{\{X_i\not\in
  (-\infty,a_kx_*)^d\}}\Big)^{2+\delta} = O(r_nv_n),
 \end{equation}
  for some $\delta>0$.
  (see Section \ref{proofs}).

  By the same arguments one can show that condition (D6) is
  fulfilled for the family $\CC:=\big\{ \big(\R^d\setminus(-\infty,x_A],\R^d\setminus(-\infty,x_B]\big),
   \big(\R^d\setminus(-\infty,x_A],\R^d\big) \mid  x_A,x_B \in (x_*,\infty)^d\big\}$.
\end{example}

From Theorem \ref{theo:extremo1} one may easily conclude the uniform
asymptotic normality of the empirical extremogram centered at the
{\em pre-asymptotic extremogram}
$$ \rho_{t,A,B}(h) := P(X_h/t\in B\mid X_0/t\in A). $$
\begin{corollary} \label{corol:empextremo}
  Suppose that the conditions of Theorem \ref{theo:extremo1} are met, that
  $(A,\R^d)\in\CC$ for all $(A,B)\in\CC$ and
  $\inf_{(A,B)\in\CC}\nu_0(A)>0$, and that $\sup_{h\in\{0,\ldots,h_0\}, (A,B)\in\CC}|\rho_{a_k,A,B}(h)-\rho_{A,B}(h)|\to 0$. Then
  \begin{eqnarray}
    \lefteqn{\sqrt{nv_n} \Big(\hat\rho_{n,A,B}(h)- \rho_{a_k,A,B}(h)\Big)_{h\in\{0,\ldots,h_0\}, (A,B)\in\CC}} \nonumber\\
    & \to &
    \Big(\frac{\nu_0\big(\R^d\setminus(-\infty,x_*)^d\big)}{\nu_0(A)}\big(\tilde
  Z(h,A,B)-\rho_{A,B}(h)\tilde Z(h,A,\R^d)\big)\Big)_{h\in\{0,\ldots,h_0\},
  (A,B)\in\CC}  \nonumber\\
    & =: & (R(h,A,B))_{h\in\{0,\ldots,h_0\},  (A,B)\in\CC}
    \label{eq:empextremoconv}
  \end{eqnarray}
  weakly.
  Hence if, in addition,
  \begin{equation}  \label{eq:preasympcond1}
    \sup_{h\in\{0,\ldots,h_0\}, (A,B)\in\CC}
    \sqrt{\frac nk}|\rho_{a_k,A,B}(h)-\rho_{A,B}(h)|\to 0,
  \end{equation}
  then
   $$
    \sqrt{nv_n} \Big(\hat\rho_{n,A,B}(h)- \rho_{A,B}(h)\Big)_{h\in\{0,\ldots,h_0\},
    (A,B)\in\CC} \to
  (R(h,A,B))_{h\in\{0,\ldots,h_0\},  (A,B)\in\CC}
  $$
  weakly.
\end{corollary}

We have already mentioned in the introduction that the empirical
extremogram $\hat\rho_{n,A,B}(h)$ is not a valid estimator if the
normalizing constants $a_k$ are unknown. In this case we replace
them by some estimator $\hat a_k$ which is consistent in the sense
that $\hat a_k/a_k\to 1$ in probability. Noting that
$$ \hat{\hat \rho}_{n,A,B}(h) := \frac{\textstyle \sum_{i=1}^n
1_{A\times B}(X_i/\hat a_k,X_{i+h}/\hat a_k)}{\textstyle
\sum_{i=1}^n 1_{A}(X_i/\hat a_k)}= \hat \rho_{n,(\hat
a_k/a_k)A,(\hat a_k/a_k)B}(h),
$$
the asymptotic normality of $\hat{\hat \rho}_{n,A,B}(h)$ is an easy
consequence of Corollary \ref{corol:empextremo}, provided that
$\rho_{t,A,B}(h)$ is a sufficiently regular function of $t$.
\begin{corollary} \label{corol:empextremo2}
  Assume that the conditions of Corollary \ref{corol:empextremo}
  (except \eqref{eq:preasympcond1}) are fulfilled and, in addition,
  $\hat a_k/a_k\to 1$ in probability, that $(A,B)\in\CC$ implies $(\lambda A,\lambda B)\in\CC$
  for all $\lambda$ in a neighborhood of 1  and that $\sup_{(A,B)\in\CC}
  \bar\varrho\big((A,B),(\lambda A,\lambda B)\big) \to 0$ as
  $\lambda\to 1$. Then
  $$
    \sqrt{nv_n} \Big(\hat{\hat \rho}_{n,A,B}(h)- \rho_{\hat a_k,A,B}(h)\Big)_{h\in\{0,\ldots,h_0\},
    (A,B)\in\CC} \to
   (R(h,A,B))_{h\in\{0,\ldots,h_0\},  (A,B)\in\CC}
  $$
  weakly. Hence, if the following second order condition holds
  \begin{equation} \label{eq:secordrho}
    \rho_{t,A,B}(h) = \rho_{A,B}(h) + \Phi_h(t)\Psi_h(A,B) +
    o(\Phi_h(t))
  \end{equation}
  uniformly for $h\in\{0,\ldots,h_0\}$, $(A,B)\in\CC$, and some extended
  regularly varying function $\Phi_h$ (see Bingham et al., 1987,
  Section 2.0) satisfying $\Phi_h(t)\to 0$ as $t\to \infty$ and some
  functions $\Psi_h$ such that $\sup_{(A,B)\in\CC}
  |\Psi_h(A,B)|<\infty$, then
   $$
     \sqrt{nv_n} \Big(\hat{\hat \rho}_{n,A,B}(h)- \rho_{a_k,A,B}(h)\Big)_{h\in\{0,\ldots,h_0\},
     (A,B)\in\CC} \to
     (R(h,A,B))_{h\in\{0,\ldots,h_0\},  (A,B)\in\CC}
     $$
  weakly, provided $\Phi_h(a_k)=O((k/n)^{1/2})$. If $\Phi_h(a_k)=o((k/n)^{1/2})$, then this
  convergence holds with $\rho_{A,B}(h)$ instead of
  $\rho_{a_k,A,B}(h)$.
\end{corollary}

\begin{remark} \label{rem:secord}
  \begin{enumerate}
    \item If $(X_0,X_h)$ satisfies the second order condition
    \begin{equation}  \label{eq:regvarsecord}
       a^\leftarrow(t) P\{(X_0,X_h)/t\in A\times B\} =
       \nu_{(0,h)}(A\times B)+\Phi_h(t)\tilde\Psi_h(A\times B)+
       o(\Phi_h(t))
    \end{equation}
    uniformly for all $(A,B)\in\CC$ with $\sup_{(A,B)\in\CC}|\tilde \Psi_h(A\times B)|<\infty$, then under the conditions of
    Corollary \ref{corol:empextremo2} direct calculations show that
    $\rho_{t,A,B}(h)=P\{(X_0,X_h)/t\in A\times B\}/P\{(X_0,X_h)/t\in A\times
    \R^d\}$ satisfies condition \eqref{eq:secordrho} with
    $\Psi_h(A, B) = \big(\tilde \Psi_h(A\times
    B)-\rho_{A,B}(h)\tilde \Psi_h(A\times\R^d)\big)/\nu_0(A).$

    \item If the conditions of Theorem \ref{theo:extremo1} hold
    when \eqref{eq:clustermoment} is replaced with
     $$ E\Big(\sum_{i=1}^{r_n} 1_{[-y,y]^d\setminus [-x,x]^d}(X_i/a_k)\Big)^2 \le u(y-x) r_n v_n
  $$
  for all $1-\delta\le x<y\le 1+\delta$, $n\in\N$, and some function $u$
  satisfying $u(t)\to 0$ as $t\downarrow 0$, then the same arguments
  as used in the proof of Theorem \ref{theo:extremo1} show that
  $$ \frac 1{\sqrt{n/k}} \sum_{i=1}^n \Big( \Ind{\|X_i\|/a_k>x}-
  P\{\|X_i\|/a_k>x\}\Big) = \sqrt{kv_n} \tilde
  Z_n(0,\R^d\setminus[-x,x]^d,\R^d), \quad x\in[1-\delta,1+\delta],
  $$
  converges weakly to a continuous Gaussian process. From
  $P\{\|X_0\|/a_k>x\}\sim x^{-\alpha}/k$ and
  \begin{eqnarray*}
    \lefteqn{ \frac{\|X\|_{n-\floor{n/k}:n}}{a_k}}\\
     & = & \inf\Big\{ x\,\Big|\,
    \sum_{i=1}^n \Ind{\|X_i\|/a_k>x}\le \floor{n/k}\Big\}  \\
    & = & \inf\Big\{ x\,\Big|\,
    \sqrt{kv_n} \tilde  Z_n(0,\R^d\setminus[-x,x]^d,\R^d)\le (k/n)^{1/2}\big(\floor{n/k}-nP\{\|X_0\|/a_k>x\}\big)
    \Big\},
  \end{eqnarray*}
  one can easily conclude that $\|X\|_{n-\floor{n/k}:n}/a_k\to 1$ in
  probability, i.e.\ $\|X\|_{n-\floor{n/k}:n}$ is consistent for
  $a_k$.

  Indeed, a refined analysis shows that under the second order
  condition \eqref{eq:regvarsecord} one even has
  $\|X\|_{n-\floor{n/k}:n}/a_k-1=O_P((k/n)^{1/2})$ if
  $\Phi(a_k)=O_P((k/n)^{1/2})$.
 \end{enumerate}
\end{remark}

As the distribution of the limit process  arising in Corollary 3.6
is difficult to estimate, we use the bootstrap approach discussed in
Section \ref{sect:bootstrap} to approximate the distribution of the
empirical extremogram. Let
\begin{eqnarray*}
  \hat\rho^*_{n,A,B}(h) & := & \frac{\sum_{j=1}^{m_n} (1+\xi_j)
  \sum_{i=1}^{r_n} 1_{A\times B} \big( a_k^{-1}
  (X_{(j-1)r_n+i},X_{(j-1)r_n+i+h})\big)}{\sum_{j=1}^{m_n} (1+\xi_j)
  \sum_{i=1}^{r_n} 1_{A} \big( a_k^{-1}
  X_{(j-1)r_n+i}\big)}\\
   \hat{\hat\rho}^*_{n,A,B}(h) & := & \frac{\sum_{j=1}^{m_n} (1+\xi_j)
  \sum_{i=1}^{r_n} 1_{A\times B} \big( \hat a_k^{-1}
  (X_{(j-1)r_n+i},X_{(j-1)r_n+i+h})\big)}{\sum_{j=1}^{m_n} (1+\xi_j)
  \sum_{i=1}^{r_n} 1_{A} \big( \hat a_k^{-1}
  X_{(j-1)r_n+i}\big)}\\
  & = & \hat\rho^*_{n,(\hat a_k/a_k)A,(\hat a_k/a_k)B}(h)\\
    R_{n,\xi}(h,A,B) & := & \sqrt{nv_n}\big( \hat\rho^*_{n,A,B}(h)-
    \hat\rho_{n,A,B}(h)\big)\\
   \hat R_{n,\xi}(h,A,B) & := & \sqrt{nv_n}\big( \hat{\hat\rho}^*_{n,A,B}(h)-
    \hat{\hat\rho}_{n,A,B}(h)\big) = R_{n,\xi}(h,(\hat a_k/a_k)A,(\hat
    a_k/a_k)B).
\end{eqnarray*}
\begin{theorem}  \label{theo:extremobootstrap}
  Suppose that all conditions of Corollary \ref{corol:empextremo2} are
  fulfilled and that $\xi_j$, $j\in\N$, are i.i.d.\ random variables
  with $E(\xi_1)=0$ and $Var(\xi_1)=1$ independent of $(X_t)_{t\in\N_0}$. Then,
  \begin{eqnarray}  \label{eq:extremobootstrap1}
   \sup_{g\in BL_1(\ell^\infty(\{0,\ldots,h_0\}\times
  \CC))} \big| E_\xi g(R_{n,\xi})-E g(R)\big| &\to & 0\\
   \sup_{g\in BL_1(\ell^\infty(\{0,\ldots,h_0\}\times
  \CC))} \big| E_\xi g(\hat R_{n,\xi})-E g(R)\big| &\to & 0
  \quad \text{in probability.}  \label{eq:extremobootstrap2}
  \end{eqnarray}
 { }
\end{theorem}

Let $\tilde \FF:=\{0,\ldots, h_0\}\times\CC$. In view of Theorem
\ref{theo:extremobootstrap}, approximate confidence regions for the
extremogram $(\rho_{A,B}(h))_{(h,A,B)\in\tilde\FF}$ can be obtained
from Monte Carlo simulations of $\hat{\hat\rho}^*_{n,A,B}(h)$. To
this end, suppose $\DD$ is a family of subsets of
$\ell^\infty(\tilde\FF)$ such that $\sup_{D\in\DD}P\{R\in
U_\eps(\partial D)\}\to 0$ as $\eps\downarrow 0$, where $U_\eps(A)$
denotes the open $\eps$-neighborhood of a set $A$ w.r.t.\ the
supremum norm $\|\cdot\|_{\tilde\FF}$ on $\ell^\infty(\tilde\FF)$.
Then all indicator functions $1_D$, $D\in\DD$, can be uniformly well
approximated from above and from below by functions of the form
$g_{\eps,A}:=(1-d_A/\eps)^+$ with $d_A(z):= \inf_{\tilde z\in A}
\|z-\tilde z\|$. Since the functions $\eps g_{\eps,A}$ belong to
$BL_1(\ell^\infty(\tilde\FF))$, it is easily seen that
\eqref{eq:extremobootstrap1} and \eqref{eq:extremobootstrap2} imply
$\sup_{D\in\DD} \big| P_\xi\{R_{n,\xi}\in D\}-P\{R\in D\}\big|\to 0$
and $\sup_{D\in\DD} \big| P_\xi\{\hat R_{n,\xi}\in D\}-P\{R\in
D\}\big|\to 0$  as $n\to\infty$, respectively.

In particular, if for sufficiently large $n\in\N$, $D_\alpha$ is a
subset of $\ell^\infty(\tilde\FF)$ such that
\begin{equation}  \label{eq:bootconf1}
 P_\xi\big\{(\hat{\hat\rho}^*_{n,A,B}(h)-
    \hat{\hat\rho}_{n,A,B}(h))_{(h,A,B)\in\tilde\FF}\in
    D_\alpha\big\}=\alpha,
\end{equation}
then under the conditions of Corollary \ref{corol:empextremo2} with
$\Phi(a_k)=o((k/n)^{1/2})$, for sufficiently large $n$,  we have
\begin{equation}  \label{eq:bootconf2}
 P\big\{
(\hat{\hat\rho}_{n,A,B}(h)-\rho_{A,B}(h))_{(h,A,B)\in\tilde\FF} \in
D_\alpha\big\} \approx \alpha.
\end{equation}
To find such a set (or rather an approximation to it), one may
simulate $B$ independent realizations
$(\hat{\hat\rho}^{*(b)}_{n,A,B}(h))_{(h,A,B)\in\tilde\FF}$, $1\le
b\le B$, of the bootstrap version of the empirical extremogram. For
some fixed set $D\subset \ell^\infty(\tilde\FF)$ let
$D_\alpha:=\lambda_\alpha D$ with $\lambda_\alpha$ denoting the
smallest $\lambda\ge 0$ such that
$$ \frac 1B \sum_{b=1}^B 1_{\textstyle
\big\{\big(\hat{\hat\rho}^{*(b)}_{n,A,B}(h)-
    \hat{\hat\rho}_{n,A,B}(h)\big)_{(h,A,B)\in\tilde\FF}\in\lambda
D\big\}}\ge\alpha.
$$
Here $D$ ought to be star-shaped, i.e.\ $z\in D$ implies $\lambda
z\in D$ for all $\lambda\in[0,1]$. The shape of $D$ determines the
emphasis which is laid on particular features of the extremogram.
See Section \ref{sect:simus} for an example.

\section{Finite sample performance of bootstrapped extremograms}
 \label{sect:simus}

In this section we investigate the finite sample performance of confidence intervals which are constructed using the multiplier block bootstrap approach, the stationary bootstrap proposed by Davis et al.\ (2012) and a modified version of the latter.

Davis et al.\ (2012) suggested to construct bootstrap samples from an observed time series $(X_t)_{1\le t\le n}$ as follows. Let $K_i$, $1\le i\le n$, be iid random variables uniformly distributed on $\{1,\ldots,n\}$, and $L_i$, $1\le i\le n$, iid random block lengths with a geometric distribution with expectation $r$, independent of $(K_i)_{1\le i\le n}$. Define $S_j:=\sum_{i=1}^j L_i$, $0\le j\le n$, $N:=\min\{j| S_j \ge n\}$, and $L_j^*:=L_j$ for $1\le j<N$ and $L_N^*:=n-S_{N-1}$. For $i\in \{S_{j-1}+1,S_{j-1}+2,\ldots, S_j\}$, $1\le j\le N$, let $X_i^*:= X_{K_j-1+i-S_{j-1}}$, where $X_t$ for $t>n$ is interpreted as $X_{(t \text{ mod } (n-1))+1}$. This means that blocks of length $L_j$ starting from the observation at $K_j$ are glued together until one obtains a new time series $(X_t^*)_{1\le t\le n}$ of length $n$; in this process one repeats the original time series after the last observation as often as necessary.
Now denote by $\hat\rho^{(*DMC)}_{n,A,B}(h)$ the bootstrap estimator of $\rho_{A,B}(h) $ calculated from $(X_t^*)_{1\le t\le n}$.
Davis et al.\ (2012) proved that under suitable conditions, conditionally on the data, the limit distribution of $\hat\rho^{(*DMC)}_{n,A,B}(h)-\hat\rho_{n,A,B}(h)$ is the same as the one of $\hat\rho_{n,A,B}(h)-\rho_{a_k,A,B}(h)$, so that bootstrap confidence intervals can be constructed.

One disadvantage of this approach is that for indices $i$ near the end of a block such that $i\le S_j< i+h$ for some $1\le j<N$ the indicator $\Ind{X^*_i\in a_kA, X^*_{i+h}\in a_kB}$ has a completely different behavior than $\Ind{X_i\in a_kA, X_{i+h}\in a_kB}$, because $(X^*_i, X^*_{i+h})$ does not correspond to a pair of observations with lag $h$ in the original time series.

To overcome this drawback, we suggest the following modification of the stationary bootstrap estimator. For simplicity, we assume that the time series has been observed at $n+h$ time points (in other words, we redefine $n$). Then we define
$$ \hat\rho_{n,A,B}^{(*stat)} (h):= \frac{\sum_{j=1}^N \sum_{i=1}^{L_j^*} \Ind{X_{K_j-1+i}\in a_k A,
  X_{K_j-1+i+h}\in a_kB}}{\sum_{j=1}^N \sum_{i=1}^{L_j^*} \Ind{X_{K_j-1+i}\in a_k A}}
$$
which has the same asymptotic behavior as $\hat\rho^{(*DMC)}_{n,A,B}(h)$, but only observations are compared which are lagged by $h$. In essence, this mean that we apply the stationary bootstrap technique to the bivariate time series $(X_t,X_{t+h})_{1\le t\le n}$.

In addition to these two version of stationary bootstrap estimators, we consider the multiplier bootstrap. Here we have drawn  multipliers  $\xi_j$ from a Student $t$-distribution with 5 degrees of freedom and scale parameter such that $Var(\xi_j)=1$. However, this particular choice is not crucial as in further simulations we have obtained a similar performance of the multiplier block bootstrap for other distributions which are symmetric about 1 with an unbounded support (e.g., for normally distributed multiplier).

Here we report the results for three different models:
\begin{enumerate}
  \item a GARCH model $X_t=\sigma_t\eps_t$, $\sigma_t^2=\alpha_0+\alpha_1X_{t-1}^2+\beta_1\sigma_{t-1}^2$ with $\alpha_0=10^{-4}, \alpha_1=0.08, \beta_1=0.9$ and $t$-distributed innovations $\eps_t$ with 8 degrees of freedom, independent of $\sigma_t$
  \item an autoregressive model of order 1: $X_t=\varphi X_{t-1}+\eps_t$ with $\varphi=0.8$ and symmetrized Fr\'{e}chet distribution of the innovations, i.e., $P\{\eps_t>x\}=P\{\eps_t<-x\}=(1-\exp(-x^{-3}))/2$ for all $x>0$
  \item a moving average time series of order 3, namely $X_t=\eps_{t}+0.5\eps_{t-1}+0.8\eps_{t-2}$ with $\eps_t$ as in (ii).
\end{enumerate}
For each model  we simulated $10\,000$ time series of length $n=2000$.

\begin{figure}[t]


\centerline{\includegraphics[width=130mm,height=180mm,
]{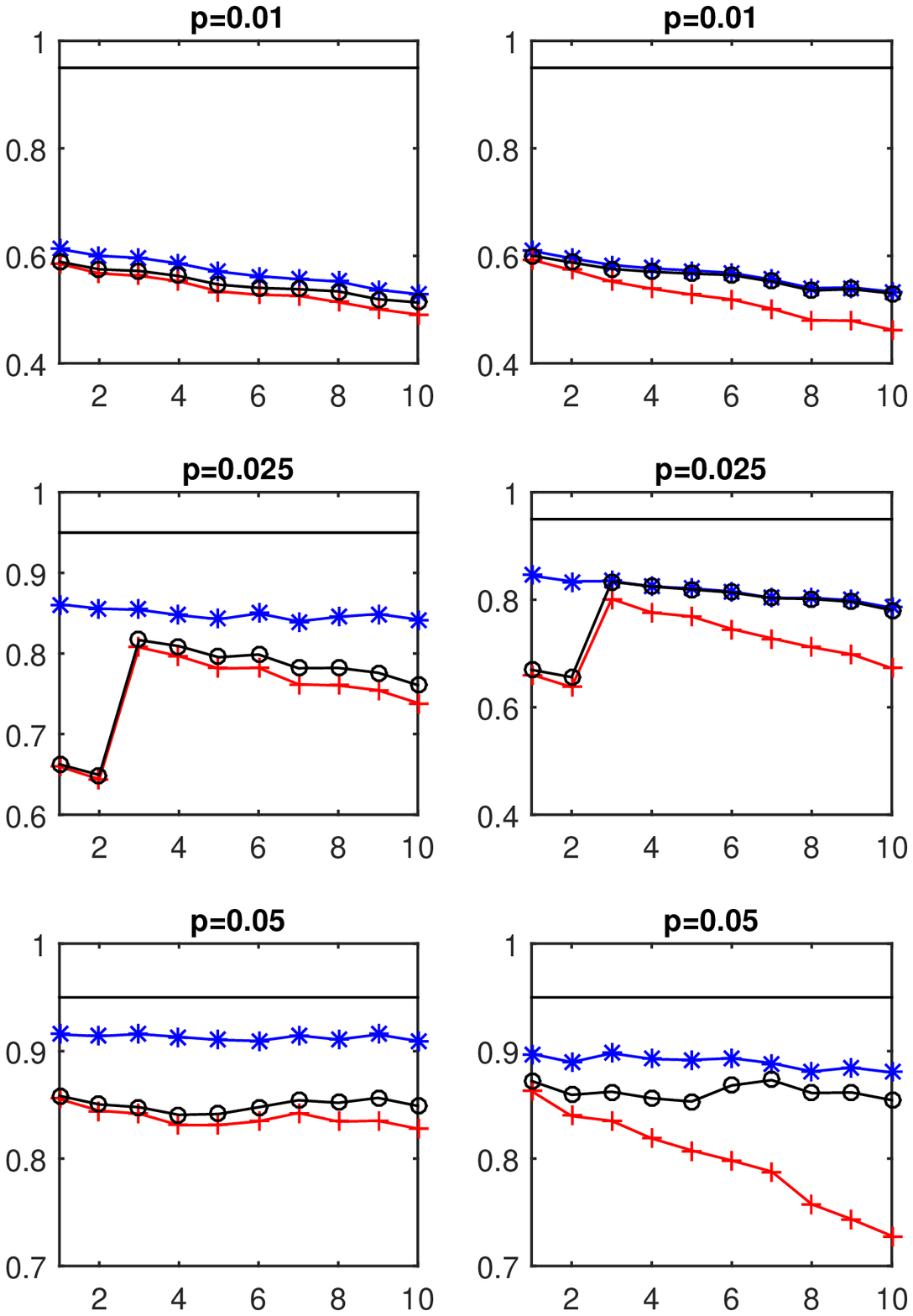}}

\caption{Empirical coverage probability of confidence intervals \eqref{eq:confint} for $\rho_{a_k,(1,\infty),(1,\infty)}(h)$ as a function of $h$, constructed using multiplier block bootstrap (blue $*$), stationary bootstrap suggested by Davis et al.\ (red $+$) and the modification thereof (black $\circ$) with (average) block length $r=100$ (left) and $r=20$ (right) for the $t$-GARCH model (i), different thresholds $a_k$ with exceedance probability $p$ are used in the three rows; the nominal coverage probability 0.95 is indicated by the horizontal line.
}
 \label{fig:covprob}
\end{figure}

We consider the extremogram for $A=B=(1,\infty)$, i.e., $\lim_{u\to\infty} P(X_h>u|X_0>u)$, which is often also called tail dependence coefficient, and lags $1\le h\le 10$. As normalizing constants $a_k$ (thresholds) we have chosen the $(1-p)$-quantile of the stationary distribution for $p\in\{0.01,0.025,0.05\}$ which have been estimated by the corresponding empirical quantiles. The true pre-asymptotic extremograms have been determined by simulation (based on 1000 time series of length $10^7$). Analytic expression for the (asymptotic) extremograms are known for the linear models (ii) and (iii) (see e.g., Meinguet and Segers, 2010, Example 9.2). For the GARCH model, they were determined using a simulation algorithm suggested by Ehlert et al.\ (2015).

In each simulation we have drawn $b=1000$ bootstrap replicas according to each of the three bootstrap procedures. If, for fixed $h$, the upper and lower empirical $\alpha/2$-quantile of the resulting $b$ bootstrap estimates of the extremogram are denoted by $u_b$ and $l_b$ then, according to \eqref{eq:bootconf1} and \eqref{eq:bootconf2},
\begin{equation} \label{eq:confint}
\big[2\hat\rho_{n,A,B}(h)-u_b, 2\hat\rho_{n,A,B}(h)-l_b\big]\cap [0,1]
\end{equation}
 is a confidence interval for the (pre-asymptotic) extremogram with nominal coverage probability $1-\alpha$.

We first discuss the results for the $t$-GARCH model in detail, before we show the results for the linear time series in abbreviated form. For this model,
Figure \ref{fig:covprob} shows the empirical coverage  probabilities of all three bootstrap procedures as a function of $h$ for the pre-asymptotic extremogram. The three rows correspond to the three thresholds with ascending exceedance probabilities. The left column shows the results for (average) block length $r=100$, the right column for $r=20$. For all bootstrap procedures, the actual coverage probabilities are much smaller than the nominal value 0.95 if the threshold is chosen too high. For the estimator based on the largest 5\% of the observations and blocks of length $r=100$, the coverage probability of the multiplier block bootstrap is reasonably close to the nominal size while both versions of the stationary bootstrap have a considerably lower coverage probability.   In all simulations, the multiplier block bootstrap yields the highest coverage probability, while the  stationary bootstrap proposed by Davis et al.\ (2012) performs worst. Moreover, in most cases the performance is better for larger block sizes. In particular, the stationary bootstrap proposed by Davis et al.\ is sensitive to too small a block size, as was to be expected from the above discussion.

\begin{figure}[tb]


\centerline{\includegraphics[width=130mm,height=60mm
]{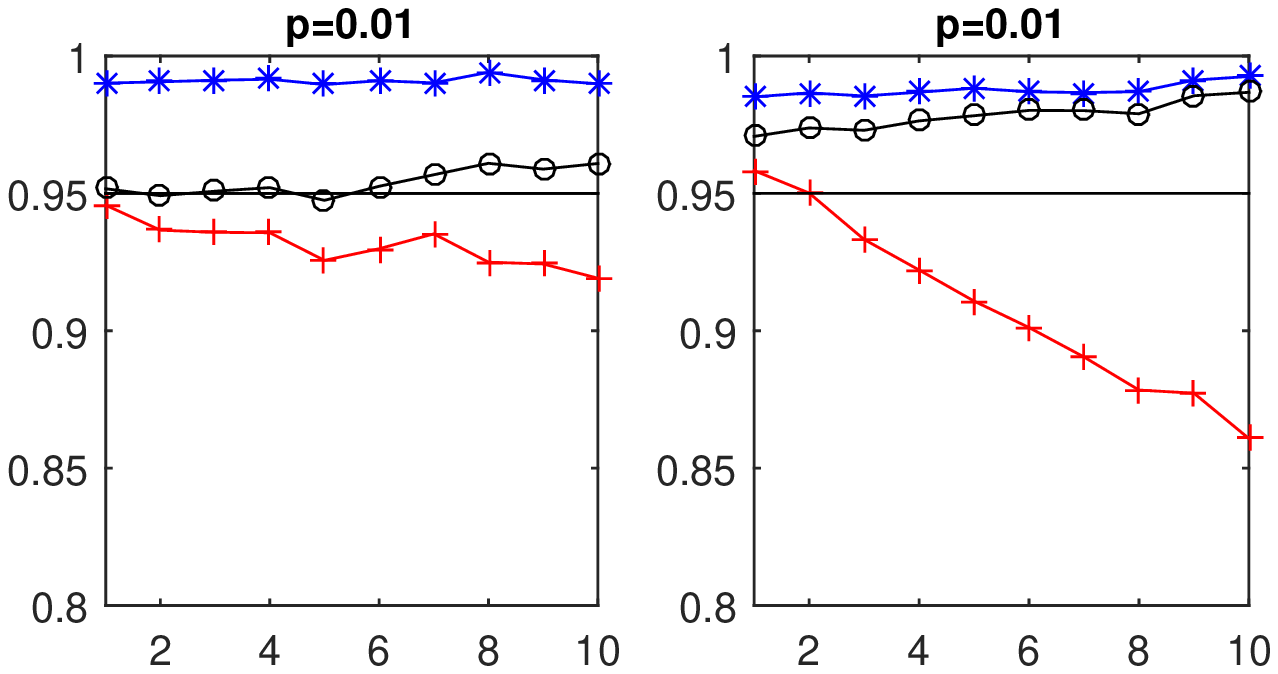}}

\caption{Empirical coverage probability of confidence intervals \eqref{eq:confint} for the pre-asymptotic extremogram $\rho_{a_k,(1,\infty),(1,\infty)}(h)$ as a function of $h$, constructed using multiplier block bootstrap (blue $*$), stationary bootstrap suggested by Davis et al.\ (red $+$) and the modification thereof (black $\circ$) with (average) block length $r=100$ (left) and $r=20$ (right) for the $t$-GARCH model (i), based only on those simulations in which for some $t$ both $X_t$ and $X_{t+h}$ exceed the threshold $a_k$.
}
 \label{fig:covprobcorr}
\end{figure}

The main reason for the disappointing performance for high thresholds is that then for very few or even none time instants both $X_t$ and $X_{t+h}$ exceed the threshold. If there are no joint exceedances in the original time series (leading to an estimate 0 for the extremogram) then also the bootstrap estimate equals 0 if one uses the multiplier block bootstrap or the modified stationary bootstrap (and it equals 0 for the original stationary bootstrap with very high probability). Hence the confidence intervals do not cover the true value if this is not exactly equal to 0, which is neither the case for  the pre-asymptotic nor the asymptotic extremogram, leading to a high non-coverage probability. Indeed, for $p=0.01$, Figure \ref{fig:covprobcorr} shows that if one considers only those simulations when the estimated extremogram does not equal 0, then the empirical coverage probability is rather close to the nominal value.

To overcome this weakness, we suggest to estimate the error distribution using a bootstrap based on a lower threshold if one wants to construct confidence intervals for the pre-asymptotic extremogram for a high threshold (or even the extremogram). Denote by $\hat\rho_{n,p}$ the empirical extremogram based on the exceedances over the threshold with exceedance probability $p$, and by $\hat\rho_{n,p}^*$ some bootstrap version thereof. Then, according to Theorem \ref{theo:extremobootstrap}, conditional on the data,  for $0<p_1<p_2$, the bootstrap error $\hat\rho_{n,p_1}^*-\hat\rho_{n,p_1}$ has approximately the same distribution as $(p_2/p_1)^{1/2}\big(\hat\rho_{n,p_2}^*-\hat\rho_{n,p_2})$. So if $u_b$ and $l_b$ denote the empirical bootstrap quantiles as defined above, calculated from the bootstrap for the threshold with the higher exceedance probability $p_2$, then
\begin{equation} \label{eq:confint2}
\Big[\hat\rho_{n,p_1}-(p_2/p_1)^{1/2}\big(u_b-\hat\rho_{n,p_2}), \hat\rho_{n,p_1}-(p_2/p_1)^{1/2}\big(l_b-\hat\rho_{n,p_2})\Big]\cap [0,1]
\end{equation}
 is a confidence interval with nominal coverage probability $1-\alpha$.

\begin{figure}[t]


\centerline{\includegraphics[width=130mm,height=120mm
]{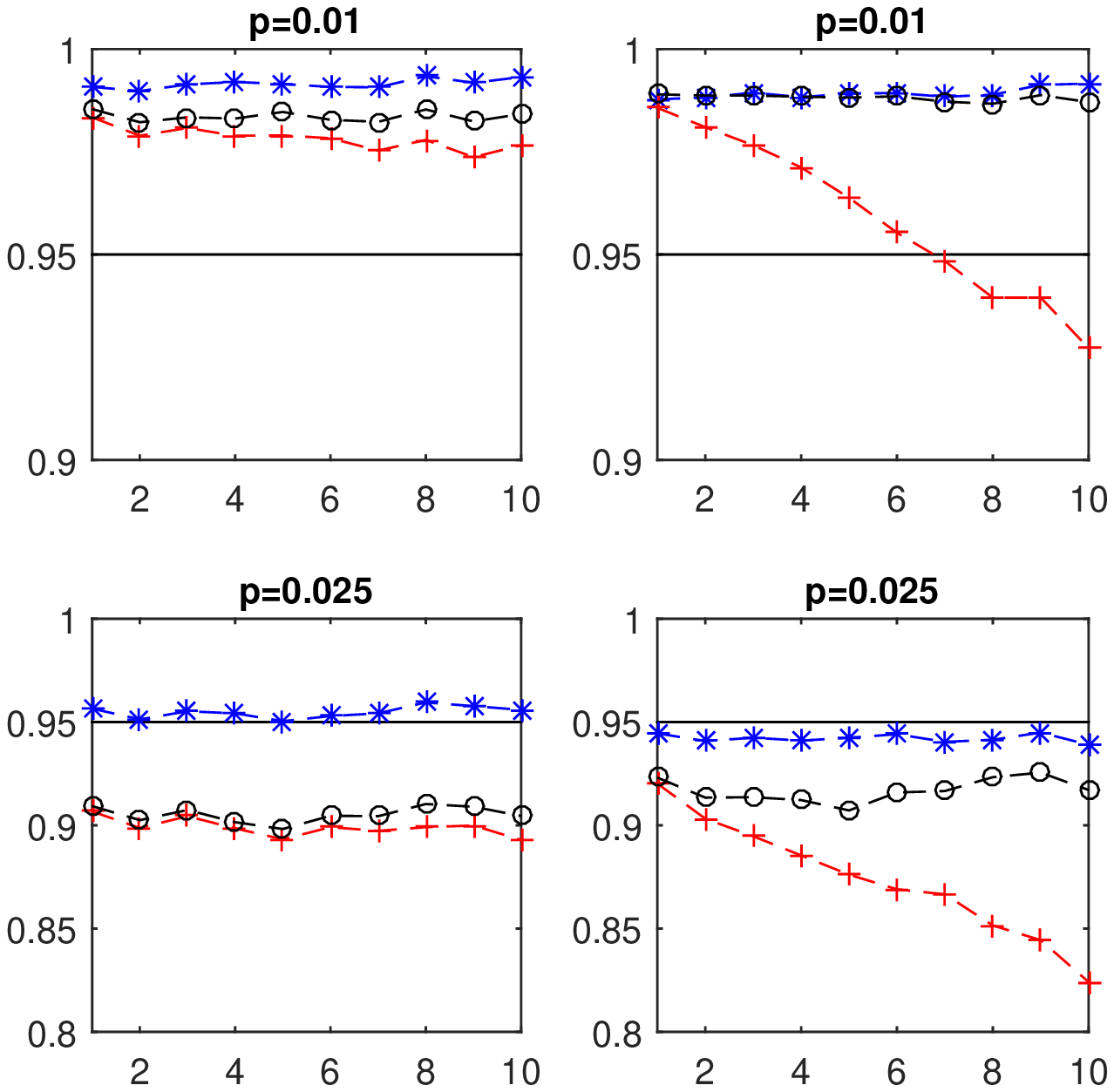}}

\caption{Empirical coverage probability of confidence intervals \eqref{eq:confint2} for $\rho_{a_k,(1,\infty),(1,\infty)}(h)$ as a function of $h$, constructed using multiplier block bootstrap (blue $*$), stationary bootstrap suggested by Davis et al.\ (red $+$) and the modification thereof (black $\circ$) with (average) block length $r=100$ (left) and $r=20$ (right) for the $t$-GARCH model (i).
}
\label{fig:covprobtrans}
\end{figure}

\begin{figure}[t]


\centerline{\includegraphics[width=130mm,height=120mm 
]{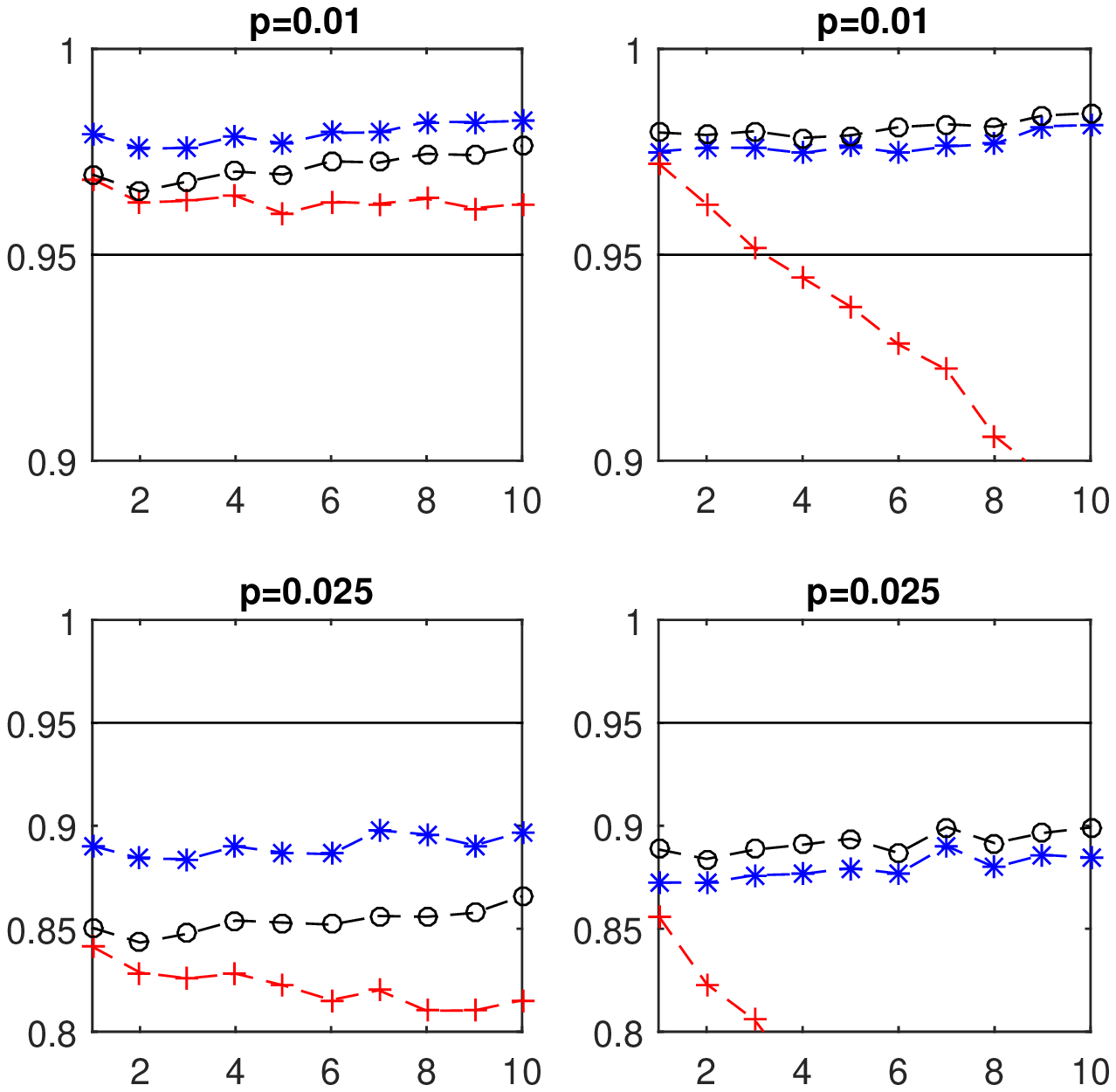}}

\caption{Empirical coverage probability of confidence intervals \eqref{eq:confint2} \eqref{eq:confint} for  the exremogram $\rho_{(1,\infty),(1,\infty)}(h)$ as a function of $h$, constructed using multiplier block bootstrap (blue $*$), stationary bootstrap suggested by Davis et al.\ (red $+$) and the modification thereof (black $\circ$) with (average) block length $r=100$ (left) and $r=20$ (right) for the $t$-GARCH model (i).
}
 \label{fig:covprobasymp}
\end{figure}

Figure \ref{fig:covprobtrans} displays the empirical coverage probabilities of this confidence interval for the pre-asymptotic extremogram, $p_1\in\{0.01,0.025\}$ and $p_2=0.05$, which are now much closer to the nominal size 0.95. (Indeed, for $p_1=0.01$ the new confidence intervals are a bit too conservative.) As for small $p_1$ the pre-asymptotic extremograms are closer to the limit extremograms, for these thresholds one may also be interested in the coverage probability for the latter, which are shown in Figure \ref{fig:covprobasymp}. The confidence intervals based on the threshold with exceedance probabilities $p_1=0.01$ are still a bit conservative, while for $p_1=0.025$, when the bias is larger and the confidence intervals more narrow, the actual coverage probabilities are too low.

\begin{figure}[t]


\centerline{\includegraphics[width=130mm,height=180mm
]{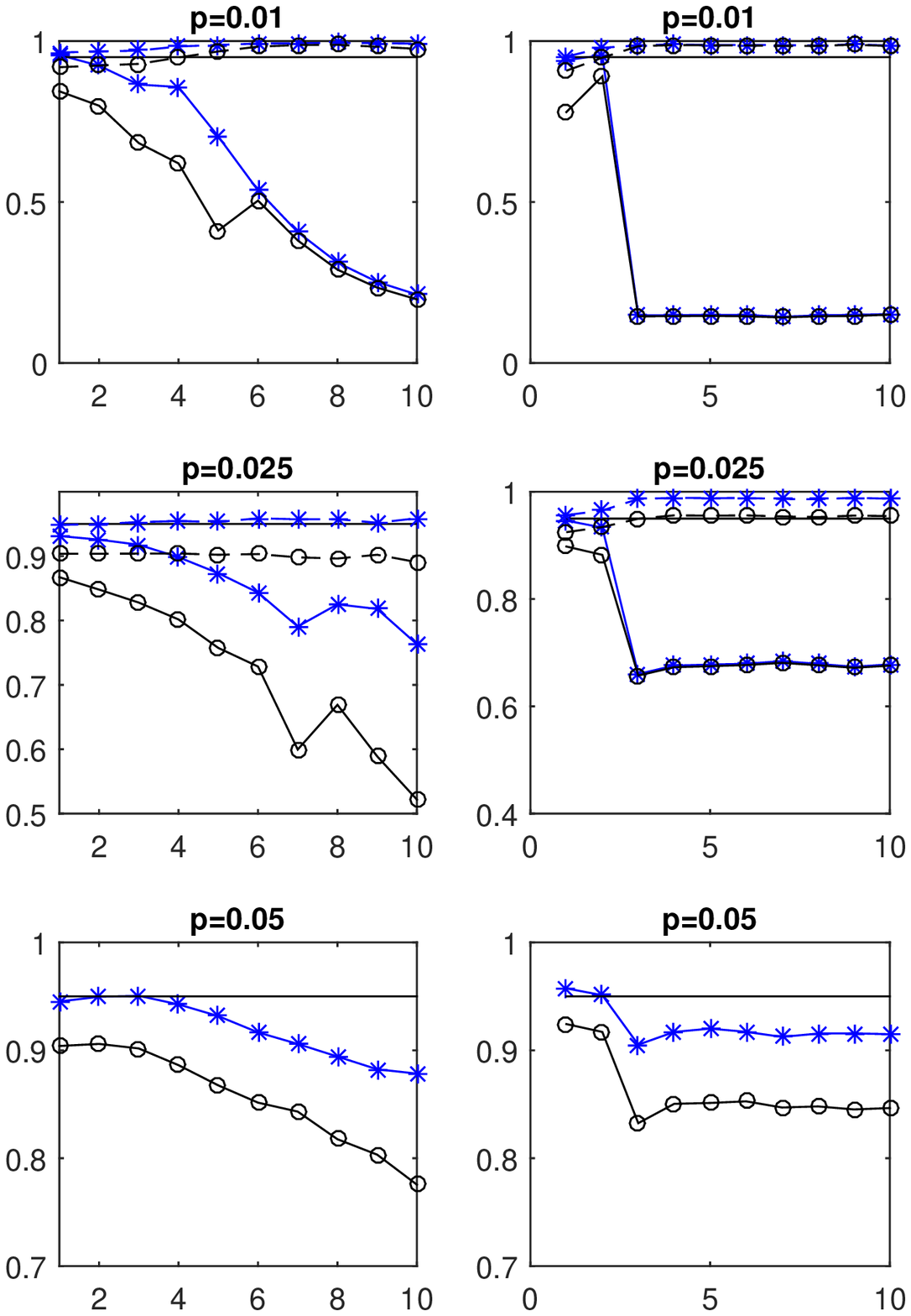}}

\caption{Empirical coverage probability of confidence intervals \eqref{eq:confint} (solid lines) and \eqref{eq:confint2} (dashed lines) constructed using multiplier block bootstrap (blue $*$) and modified stationary bootstrap (black $\circ$) with  block length $r=100$ for the AR(1) model (ii) (left) and MA(3) model (iii) (right) and different thresholds with exceedance probability $p$; the nominal coverage probability 0.95 is indicated by the horizontal line.
}
\label{fig:covproblinear}
\end{figure}

Finally, we briefly discuss the linear time series models. As overall the conclusions are similar, we present just the most important findings for block size $r=100$. Figure \ref{fig:covproblinear} shows the coverage probabilities for the autoregressive model (ii) in the left column and for the moving average (iii) in the right column, both for the confidence intervals \eqref{eq:confint} (solid lines) and \eqref{eq:confint2} (dashed lines). In order to  not overload the plot, the results for the original stationary bootstrap (which again performed worst) are not shown. Again the multiplier block bootstrap gives the highest coverage probabilities, which are nevertheless not satisfactory if one uses the direct bootstrap interval
\eqref{eq:confint} for a high threshold for the extremogram at lags not close to 0. This is particularly true, if the true value is small (e.g., for large lags in the autoregressive model). In these cases, it helps a lot to borrow strength from the bootstrap for a lower threshold as in \eqref{eq:confint2}.

\section{Proofs}
  \label{proofs}

\begin{proofof} Theorem \ref{theo:bootfidis}. \rm\quad
  We combine ideas from the proofs of Theorem 2.3 of Drees and
  Rootz\'{e}n (2010) and of Theorem 2 by Kosorok (2003). Denote by
  $Y_{n,j}^*$, $1\le j\le m_n$, independent copies of $Y_{n,j}$ that
  are independent of $(\xi_i)_{i\in\N}$. As in Drees and
  Rootz\'{e}n (2010), we define $ \Delta^*_{n,j}(f)  :=   f(Y_{n,j}^*) - f((Y_{n,j}^*)^{(r_n-l_n)})$, $1\le j\le
  m_n$. (Recall that $x^{(l)}:=(x_1,\ldots,x_l)$ for
  $x=(x_1,\ldots,x_r)$ with $r\ge l$.)

  We first analyze the asymptotic behavior of
  $$    \frac 1{\sqrt{n v_n}} \sum_{j=1}^{\floor{m_n/2}} \big( \Delta_{n,2j}^*(f)-
E \Delta_{n,2j}^*(f)\big)
  $$
  conditionally given $(Y_{n,j}^*)_{1\le j\le m_n}$. Note that $E_\xi\big(\xi_{2j} (\Delta_{n,2j}^*(f)-
E \Delta_{n,2j}^*(f))\big)=E\big( \xi_{2j} (\Delta_{n,2j}^*(f)-
E \Delta_{n,2j}^*(f))\mid (Y^*_{n,j})_{1\le j\le m_n}\big)=0$. Moreover, because of $E\xi_{2j}^2=1$
  \begin{eqnarray}
    \lefteqn{   \frac 1{n v_n} \sum_{j=1}^{\floor{m_n/2}} E_\xi \Big( \xi_{2j}^2 (\Delta_{n,2j}^*(f)-
E \Delta_{n,2j}^*(f))^2\Ind{|\xi_{2j} (\Delta_{n,2j}^*(f)-
E \Delta_{n,2j}^*(f))|\le \sqrt{nv_n}}\Big)} \nonumber \\
  & \le &  \frac 1{n v_n} \sum_{j=1}^{\floor{m_n/2}}(\Delta_{n,2j}^*(f)-E \Delta_{n,2j}^*(f))^2.\hspace*{6cm}
   \label{eq:ineq1}
  \end{eqnarray}
  Now by condition (C1)
  \begin{eqnarray}
    \lefteqn{P\Big\{ \sum_{j=1}^{\floor{m_n/2}} (\Delta_{n,2j}^*(f)-E \Delta_{n,2j}^*(f))^2
    \Ind{|\Delta_{n,2j}^*(f)-E \Delta_{n,2j}^*(f)|> \sqrt{nv_n}} \ne 0\Big\}} \nonumber\\
     & \le & \floor{m_n/2} P\{ |\Delta_{n}^*(f)-E \Delta_{n}^*(f)| > \sqrt{nv_n}\} \nonumber\\
     & \to & 0 \hspace*{10cm}\label{eq:conv1}
   \end{eqnarray}
  and
  \begin{eqnarray}
   \lefteqn{E\Big( \frac 1{n v_n} \sum_{j=1}^{\floor{m_n/2}}(\Delta_{n,2j}^*(f)-E \Delta_{n,2j}^*(f))^2
  \Ind{|\Delta_{n,2j}^*(f)-E \Delta_{n,2j}^*(f)|\le \sqrt{nv_n}}\Big) }\nonumber \\
    & \le & \frac {m_n}{2n v_n} E\Big((\Delta_{n}^*(f)-E \Delta_{n}^*(f))^2
  \Ind{|\Delta_{n}^*(f)-E \Delta_{n}^*(f)|\le \sqrt{nv_n}}\Big) \nonumber \\
  & = & o\Big(\frac{r_n v_n m_n}{n v_n} \Big) \nonumber\\
  & = & o(1).  \label{eq:conv2}
  \end{eqnarray}
  Combining \eqref{eq:ineq1}--\eqref{eq:conv2}, we see that the
  left-hand side of \eqref{eq:ineq1} tends to 0 in probability.

  Next check that from \eqref{eq:conv1} and \eqref{eq:conv2}  we may conclude
  \begin{eqnarray}
    \lefteqn{\sum_{j=1}^{\floor{m_n/2}} P_\xi \big\{|\xi_{2j} (\Delta_{n,2j}^*(f)-
E \Delta_{n,2j}^*(f))|> \sqrt{nv_n}\big\} } \nonumber\\
   & \le & \sum_{j=1}^{\floor{m_n/2}}  \Big( \frac{E(\xi_{2j}^2) (\Delta_{n,2j}^*(f)-
E \Delta_{n,2j}^*(f))^2}{nv_n} \Ind{|\Delta_{n,2j}^*(f)-E \Delta_{n,2j}^*(f)|\le \sqrt{nv_n}}  \nonumber\\
  & & \hspace*{2cm} + \Ind{|\Delta_{n,2j}^*(f)-E \Delta_{n,2j}^*(f)|> \sqrt{nv_n}}\Big) \nonumber\\
  & \to & 0 \label{eq:conv3}
  \end{eqnarray}
  in probability. Therefore, to each subsequence $n'$ there exists a subsubsequence
   $n''$ such that the convergence of the left-hand side of \eqref{eq:ineq1} and the convergence
   of the left-hand side of
   \eqref{eq:conv3} hold almost surely. By Theorem 4.10 of Petrov (1995), on the corresponding set of probability
   1, for all $\eta>0$,
   \begin{equation} \label{eq:conv4}
      P_\xi\Big\{ \frac 1{\sqrt{n'' v_{n''}}} \Big|\sum_{j=1}^{\floor{m_{n''}/2}} \xi_{2j} \big( \Delta_{n'',2j}^*(f)-
E \Delta_{n'',2j}^*(f)\big) \Big| >\eta\Big\} \;\longrightarrow\; 0.
   \end{equation}

  We can argue the same way to obtain convergence \eqref{eq:conv4}
  uniformly for a finite number of cluster functionals $f_1,\ldots, f_l$ and
  for the analogous sum over the odd numbered blocks.

  By Lemma 3 of Kosorok (2003) and the conditions (C2) and (C3), the subsubsequence $n''$ can be
  chosen such that on a set with probability 1
  $$ \sup_{g\in BL_1(\R^l)} \Big| E_\xi g\Big(\frac 1{\sqrt{n'' v_{n''}}}
  \sum_{j=1}^{m_{n''}} \xi_j \big( f_k(Y_{n'',j}^*)- E f_k(Y_{n'',j}^*)\big)_{1\le k\le
  l}\Big)
   - Eg(Z((f_k)_{1\le k\le l}))\Big| \to 0.
  $$
  Because of the aforementioned generalizations of \eqref{eq:conv4} it follows that
  \begin{eqnarray*}
    \lefteqn{ \sup_{g\in BL_1(\R^l)} \Big| E_\xi g\Big(\frac 1{\sqrt{n'' v_{n''}}}
  \sum_{j=1}^{m_{n''}} \xi_j \big( f_k((Y_{n'',j}^*)^{(r_{n''}-l_{n''})})-
  E f_k((Y_{n'',j}^*)^{(r_{n''}-l_{n''})})\big)_{1\le k\le l}}\\
   & & \hspace*{8cm}
   - Eg\big(Z((f_k)_{1\le k\le l})\big)\Big| \to 0. \hspace*{1cm}
  \end{eqnarray*}
  Since, by (B2),
  $$  \big\| P^{(Y_{n,j}^{(r_n-\ell_n)})_{1\le j\le
    m_n}} - P^{((Y^*_{n,j})^{(r_n-\ell_n)})_{1\le j\le
    m_n}} \big\|_{TV} \le m_n \beta_{n,l_n}
    \; \longrightarrow\; 0
  $$
  (see Drees and Rootz\'{e}n, 2010, proof of Lemma 5.1), the last convergence in turn implies
  \begin{equation} \label{eq:conv5}
     \sup_{g\in BL_1(\R^l)} \Big| E_\xi g\Big(\frac 1{\sqrt{n'' v_{n''}}}
  \sum_{j=1}^{m_{n''}} \xi_j \big( f_k(Y_{n'',j}^{(r_{n''}-l_{n''})})
  - E f_k(Y_{n'',j}^{(r_{n''}-l_{n''})})\big)_{1\le k\le l}\Big)
   - Eg(Z((f_k)_{1\le k\le l}))\Big| \; \longrightarrow\;0
  \end{equation}
  in probability. Hence, along a further subsequence of $n''$, the convergence
  holds almost surely and w.l.o.g.\ we may
  assume almost sure convergence along $n''$.

  By the above arguments, one easily sees that the analog to \eqref{eq:conv4}
  also holds for $\Delta_{n,2j}(f)$ instead of $\Delta_{n,2j}^*(f)$. Together with the same
  argument for the odd numbered blocks it follows that $n''$ can be chosen such that on a set with probability
  1, for all $\eta>0$,
  $$ P_\xi\Big\{ \Big\| \Big(\frac 1{\sqrt{n'' v_{n''}}} \sum_{j=1}^{m_{n''}} \xi_j
  \big( \Delta_{n'',j}(f_k)-
E \Delta_{n'',j}(f_k)\big) \Big)_{1\le k\le l}\Big\| >\eta\Big\} \;
\longrightarrow\; 0.
  $$
  Thus from \eqref{eq:conv5} we can conclude that for all subsequences $n'$
  there exists a subsubsequence $n''$ such that almost surely
  $$ \sup_{g\in BL_1(\R^l)} \Big| E_\xi g\Big(\frac 1{\sqrt{n'' v_{n''}}}
  \sum_{j=1}^{m_{n''}} \xi_j \big( f_k(Y_{n'',j})- E f_k(Y_{n'',j})\big)_{1\le k\le l}\Big)
   - Eg\big(Z((f_k)_{1\le k\le l})\big)\Big| \; \longrightarrow\; 0,
  $$
  which is equivalent to the assertion.
\end{proofof}

\begin{proofof} Proposition  \ref{prop:boottight}. \quad \rm
  The asymptotic tightness of $Z_{n,\xi}$ follows if we can prove
  asymptotic tightness of $\big((nv_n)^{-1/2}
  \sum_{j=1}^{\floor{m_n/2}} \xi_{2j}(f(Y_{n,2j})-
  Ef(Y_{n,2j}))\big)_{f\in\FF}$ and the analogous assertion for the
  sum over the odd numbered blocks. Similarly as in the proof of
  Theorem 2.8 of Drees and Rootz\'{e}n (2010), it suffices to prove
  tightness of $\big((nv_n)^{-1/2}
  \sum_{j=1}^{\tilde m_n} \xi_{j}(f(Y_{n,j}^*)-
  Ef(Y_{n,j}^*))\big)_{f\in\FF}$ with $\tilde
  m_n\in\{\floor{m_n/2},\ceil{m_n/2}\}$ and $Y_{n,j}^*$ denoting
  independent copies of $Y_{n,j}$, because the total variation
  distance between the distribution of the processes with dependent
  blocks (which are separated in time)
  resp.\ with independent blocks tends to 0.
To this end, we verify that the conditions of van der Vaart and
Wellner (1996), Theorem 2.11.9, are fulfilled for
$Z_{ni}:=(nv_n)^{-1/2}\xi_i (f(Y_{n,i}^*)-  Ef(Y_{n,i}^*))$ which
are centered random variables because of the independence of $\xi_i$
and $Y_{n,i}^*$.

The second displayed formula of this theorem is an immediate
consequence of condition (D3), since $E \xi_i^2=1$ implies
$E\big(\xi_i(f(Y_{n,i}^*)- g(Y_{n,i}^*))\big)^2=E\big(f(Y_{n,i}^*)-
g(Y_{n,i}^*)\big)^2$. Likewise, the bracketing number for the
multiplier process considered here is the same as the bracketing
number for the original process so that the bracketing entropy
condition (i.e.\ the third displayed formula in Theorem 2.11.9)
follows from (D4).

It remains to verify that
\begin{equation}  \label{eq:L1Lindeberg}
 \frac{\tilde m_n}{\sqrt{nv_n}} E^*\big| \xi F(Y_n) \Ind{|\xi
F(Y_n)|>\eta\sqrt{nv_n}}\big| \;\longrightarrow\; 0,\quad \forall\,
\eta>0.
\end{equation}
If $\xi$ is bounded, then this convergence is obvious from (D2).

Under the conditions of part (ii), one has for all $u_n>0$
\begin{eqnarray*}
  \lefteqn{E^* \big| \xi^2 F^2(Y_n) \Ind{|\xi
F(Y_n)|>\eta\sqrt{nv_n}}\big|} \\
  & \le & E\big(\xi^2 \Ind{|\xi|>u_n}\big) E^*(F^2(Y_n)) + E\big(\xi^2 \Ind{|\xi|\le
  u_n}\big)E^*\Big(F^2(Y_n)\Ind{F(Y_n)>\eta\sqrt{nv_n}/u_n}\Big).
\end{eqnarray*}
By condition (D2') one can find a sequence $u_n\to\infty$ such that
$$ E\Big(F^2(Y_n)\Ind{F(Y_n)>\eta\sqrt{nv_n}/u_n}\Big) =o(r_nv_n). $$
Moreover, also the first term is of smaller order than $r_nv_n$,
because $E\big(\xi^2 \Ind{|\xi|>u_n}\big)\to 0$ and, by assumption,
$E^*(F^2(Y_n))=O(r_nv_n)$. Now, by the Cauchy-Schwarz inequality and
the Chebyshev inequality, the left-hand side of
\eqref{eq:L1Lindeberg} can be bounded by
\begin{eqnarray*}
  \lefteqn{ \frac{\tilde m_n}{\sqrt{nv_n}}\Big( E^*\big| \xi^2 F^2(Y_n)
\Ind{|\xi F(Y_n)|>\eta\sqrt{nv_n}}\big|\cdot E^*\Ind{|\xi
F(Y_n)|>\eta\sqrt{nv_n}}\Big)^{1/2}}\\
 & \le & o\Big( \frac{\tilde
m_n}{\sqrt{nv_n}} (r_nv_n)^{1/2}\Big)\Big(\frac{E^*(\xi^2
F^2(Y_n))}{\eta^2 nv_n}\Big)^{1/2}
 \;\longrightarrow\; 0. \hspace*{4cm}
\end{eqnarray*}
\end{proofof}

\begin{proofof} Theorem \ref{theo:bootconv1}. \quad \rm
By (D3) the family $\FF$ is totally bounded w.r.t.\ the metric
$\rho$. Hence there exists a sequence of finite $\delta$-nets
$\FF_\delta$ of $\FF$, i.e.\ finite sets such that to every
$f\in\FF$ there exists $\pi_\delta(f)\in\FF_\delta$ whose
$\rho$-distance to $f$ is less than $\delta$. Because $Z$ has
continuous sample paths w.r.t.\ $\rho$ and $g\in
BL_1(\ell^\infty(\FF))$ is bounded and Lipschitz-continuous with
Lipschitz-constant 1, we may conclude
\begin{equation}  \label{eq:gZdiscret}
  \lim_{\delta\downarrow 0} E^* \sup_{g\in BL_1(\ell^\infty(\FF))}
  \big| g(Z(\pi_\delta\circ \cdot))-g(Z(\cdot))\big| = 0.
\end{equation}

For fixed $\delta>0$, denote by $l=\sharp \FF_\delta$ the
cardinality of the $\delta$-net. Theorem \ref{theo:bootfidis} gives
\begin{eqnarray}
  \lefteqn{ \sup_{g\in BL_1(\ell^\infty(\FF))} \big| E_\xi
  g(Z_{n,\xi}(\pi_\delta\circ \cdot))- E g(Z(\pi_\delta\circ\cdot))\big|}
  \nonumber\\
  & \le & \sup_{h\in BL_1(\R^l)} \big|E_\xi
  h\big((Z_{n,\xi}(f))_{f\in\FF_\delta}\big)- E
  h\big((Z(f))_{f\in\FF_\delta}\big)\big|\nonumber\\
  & \to & 0   \label{eq:Znxidiskretconv}
\end{eqnarray}
in outer probability (cf.\ van der Vaart and Wellner, 1996, p.\
182).

Next note that by the definition of $BL_1(\ell^\infty(\FF))$
$$
\sup_{g\in BL_1(\ell^\infty(\FF))} \big| E_\xi
  g(Z_{n,\xi}(\pi_\delta\circ \cdot))- E_\xi g(Z_{n,\xi})\big|
  \le  E_\xi \min\Big(
  \sup_{f\in\FF} |Z_{n,\xi}(\pi_\delta(f))-Z_{n,\xi}(f)|, 2\Big).
$$
Since $Z_{n,\xi}$ weakly converges to $Z$, it is asymptotically
equicontinuous, that is, for all $\eps>0$ and all sequences
$\delta_n\downarrow 0$
$$   P^* \Big\{
\sup_{f,g\in\FF, \rho(f,g)<\delta_n}
|Z_{n,\xi}(f)-Z_{n,\xi}(g)|>\eps\Big\} \;\longrightarrow\; 0.
$$
Hence
$$ E^* \min\Big(\sup_{f,g\in\FF, \rho(f,g)<\delta_n}
|Z_{n,\xi}(f)-Z_{n,\xi}(g)|,2\Big) \;\longrightarrow\; 0,
$$
 and thus by Fubini's theorem (van der Vaart and Wellner, 1996,
 Lemma 1.2.6)
$$ E^* \Big(E_\xi \min\Big(\sup_{f\in\FF}
|Z_{n,\xi}(\pi_{\delta_n}(f))-Z_{n,\xi}(f)|,2\Big)\Big)
\;\longrightarrow\; 0.
$$
This in turn implies
\begin{equation}  \label{eq:Znxidiscret}
  \sup_{g\in BL_1(\ell^\infty(\FF))} \big| E_\xi
  g(Z_{n,\xi}(\pi_{\delta_n}\circ \cdot))- E_\xi g(Z_{n,\xi})\big| \;\longrightarrow\; 0
\end{equation}
in outer probability for all $\delta_n\downarrow 0$.

By \eqref{eq:gZdiscret}, for all $\eps>0$ and all
$\delta_n\downarrow 0$, one has for sufficiently large $n$ that $E^*
\sup_{g\in BL_1(\ell^\infty(\FF))}
  \big| g(Z(\pi_{\delta_n}\circ \cdot))-g(Z(\cdot))\big|<\eps/3$. Therefore, in view of
\eqref{eq:Znxidiskretconv} and \eqref{eq:Znxidiscret}, for all
$\eps,\eta>0$ and sufficiently large $n$
\begin{eqnarray*}
 \lefteqn{P^*\Big\{\sup_{g\in BL_1(\ell^\infty(\FF))} \big| E_\xi
 g(Z_{n,\xi})- Eg(Z)| >\eps\Big\}}\\
 & \le & P^*\Big\{\sup_{g\in BL_1(\ell^\infty(\FF))} \big| E_\xi
 g(Z_{n,\xi}(\cdot))-E_\xi  g(Z_{n,\xi}(\pi_{\delta_n}\circ
 \cdot))\big|>\eps/3\Big\} \\
 & & + P^*\Big\{\sup_{g\in BL_1(\ell^\infty(\FF))} \big| E_\xi  g(Z_{n,\xi}(\pi_{\delta_n}\circ
 \cdot))- E g(Z(\pi_{\delta_n}\circ \cdot))\big| > \eps/3\Big\} \\
 & < & \eta,
\end{eqnarray*}
which proves the assertion.
\end{proofof}

\begin{proofof} Proposition \ref{prop:bootequi}. \rm\quad
  For $f\in\FF$ define $Tf:\R\times E_\cup \to \R$, $Tf(t,y):=t
  f(y)$ and $T\FF:=\{Tf\mid f\in\FF\}$. We are going to apply Theorem 2.11.1 of van der Vaart and
  Wellner (1996) to the processes
  $$ \tilde Z_n(g) :=    \frac 1{\sqrt{n v_n}} \sum_{j=1}^{\tilde m_n} \big( g(\xi_j,Y_{n,j}^*)-
E g(\xi_j,Y_{n,j}^*)\big), \quad g\in T\FF,
  $$
  with $\tilde  m_n\in\{\floor{m_n/2},\ceil{m_n/2}\}$ and $Y_{n,j}^*$ denoting
  independent copies of $Y_{n,j}$. The assertion then follows by
  the same arguments as used in the proof of Drees and Rootz\'{e}n
  (2010), Theorem 2.10 (cf.\ also the proof of Proposition
  \ref{prop:boottight} of the present paper).

  Because $|TF|$ is an envelope function of $T\FF$ and
  \begin{eqnarray*}
    \lefteqn{\sum_{j=1}^{\tilde m_n} E^*\Big( \big((nv_n)^{-1/2}
    TF(\xi_j,Y_{n,j}^*)\big)^2 \Ind{(nv_n)^{-1/2}
    |TF(\xi_j,Y_{n,j}^*)|>\eta}\Big)}\\
    & \le & \frac 1{r_nv_n} E^*\big( \xi^2 F^2(Y_n)
    \Ind{|\xi F(Y_n)|>\eta\sqrt{nv_n}}\big),
  \end{eqnarray*}
  the first condition of Theorem 2.11.1 is obviously fulfilled if
  $\xi$ is bounded and (D2') holds, while it follows from the
  arguments given at the end of the proof of Proposition
  \ref{prop:boottight} if $E(F^2(Y_n))=O(r_nv_n)$ holds.

  The second condition of Theorem 2.11.1 is equivalent to our
  condition (D3), because $E(\xi_j^2)=1$ and the independence of
  $\xi_j$ and $Y_{n,j}^*$ imply
  $E\big(Tf(\xi_j,Y_{n,j}^*)-Tg(\xi_j,Y_{n,j}^*)\big)^2=E\big(f(Y_{n,j}^*)-g(Y_{n,j}^*)\big)^2$.

  It remains to verify the metric entropy condition (2.11.2) of van der Vaart and
  Wellner (1996), which is equivalent to
  $$ \lim_{\delta\downarrow 0} \limsup_{n\to\infty} P^* \Big\{
  \int_0^\delta \sqrt{\log N(\eps,T\FF,\tilde d_n)}\, d\eps>\eta\Big\}=0
  $$
  for all $\eta>0$ where
  $$ \tilde d_n(Tf,Tg) := \Big( \frac 1{n v_n} \sum_{j=1}^{\tilde
  m_n} \xi_j^2 \big( f(Y^*_{n,j})-g(Y^*_{n,j})\big)^2\Big)^{1/2}.
  $$

  If $|\xi|\le c$, then
   $  \tilde d_n(Tf,Tg)\le c d_n(f,g)$ so that $N(\eps,T\FF,\tilde
   d_n)\le N(\eps/c,\FF,d_n)$ and the entropy condition readily
   follows from (D6).

  If $\xi$ is not necessarily bounded, but the uniform entropy
  condition (D6') holds, then one may proceed similarly as in the
  proof of Theorem 2.10 of Drees and Rootz\'{e}n (2010). Let
  $$ Q_{n,\xi} := \frac{\sum_{j=1}^{\tilde m_n} \xi_j^2
  \eps_{Y_{n,j}^*}}{\sum_{j=1}^{\tilde m_n} \xi_j^2} \in \mathcal{Q}
  $$
  with $\eps_y$ denoting the Dirac measure with mass 1 at $y$,
  and check that
  $$ \tilde d_n(Tf,Tg) = \bigg(\frac{\sum_{j=1}^{\tilde m_n}
  \xi_j^2}{nv_n}\bigg)^{1/2} d_{Q_{n,\xi}}(f,g).
  $$
  Hence $N(\eps,T\FF,\tilde d_n)\le
  N\big(\eps(nv_n/\sum_{j=1}^{\tilde m_n} \xi_j^2)^{1/2},\FF,
  d_{Q_{n,\xi}}\big)$. Moreover, for all $\tau>0$
  \begin{eqnarray*}
    P\Big\{ \Big(\int F^2\, d Q_{n,\xi}\Big)^{1/2} > \tau \Big(\frac{nv_n}{\sum_{j=1}^{\tilde m_n}
    \xi_j^2}\Big)^{1/2}\Big\}
    & = & P\Big\{ \sum_{j=1}^{\tilde m_n} \xi_j^2
    F^2(Y_{n,j}^*)>\tau^2 nv_n\Big\} \\
    & \le & \frac 1{\tau^2 nv_n} E\Big(\sum_{j=1}^{\tilde m_n} \xi_j^2
    F^2(Y_{n,j}^*)\Big)\\
    & \le & \frac{E(F^2(Y_n))}{\tau^2 r_n v_n}.
  \end{eqnarray*}
  Since $E(F^2(Y_n))=O(r_nv_n)$, this probability can be made
  arbitrarily small for all $n$ by choosing $\tau$ sufficiently
  large.
  Thus, for all $\eta>0$, there exists $\tau>0$ such that with outer
  probability of at least $1-\eta$
  \begin{eqnarray*}
    \int_0^\delta \sqrt{\log N(\eps,T\FF,\tilde d_n)}\, d\eps
    & = & \tau \int_0^{\delta/\tau} \sqrt{\log N(\eps\tau,T\FF,\tilde d_n)}\,
    d\eps\\
    & \le &  \tau \int_0^{\delta/\tau} \sqrt{\log N\Big(\eps\tau\Big( \frac{nv_n}{\sum_{j=1}^{\tilde m_n}
    \xi_j^2}\Big)^{1/2},\FF,d_{Q_{n,\xi}}\Big)}\,
    d\eps\\
   & \le &  \tau \int_0^{\delta/\tau} \sup_{Q\in\mathcal{Q}}
   \sqrt{\log N\Big(\eps\Big(\int F^2\, dQ\Big)^{1/2},\FF,
   d_Q\Big)}\,    d\eps\\
   & \to & 0
  \end{eqnarray*}
  as $\delta\downarrow 0$ by (D6').

  Hence, under both sets of conditions, the asymptotic
  equicontinuity follows from Theorem 2.11.1 of van der Vaart and
  Wellner (1996).
\end{proofof}

\begin{proofof} Corollary \ref{corol:Znxistarconv}.\rm\quad
  Because
  $$ Z_{n,\xi}^*(f)-Z_{n,\xi}(f) = \frac 1{\sqrt{nv_n}}
  \sum_{j=1}^{m_n} \xi_j(Ef(Y_{n,j})-\overline{f(Y_n)}) = -\frac
  1{m_n} \sum_{j=1}^{m_n} \xi_j \cdot Z_n(f),
  $$
    $$ E_\xi \Big|\frac  1{m_n} \sum_{j=1}^{m_n} \xi_j\Big| \le \Big(E_\xi \Big(\frac  1{m_n} \sum_{j=1}^{m_n}
  \xi_j\Big)^2\Big)^{1/2} = \Big(\frac  1{m_n} Var(\xi)\Big)^{1/2}
  = \frac 1{\sqrt{m_n}}
  $$
  and $Z_n\to Z$ weakly in $\ell^\infty(\FF)$, one has
  $$  E_\xi \sup_{f\in\FF}   |Z_{n,\xi}^*(f)-Z_{n,\xi}(f)| \le E_\xi
      \Big|\frac  1{m_n} \sum_{j=1}^{m_n} \xi_j\Big|\cdot
      \sup_{f\in\FF} |Z_n(f)| \;\longrightarrow\; 0,
  $$
  which implies \eqref{eq:ZnxiZnxistarapprox}.
  Hence the weak convergence $Z_{n,\xi}^*\to Z$ follows from the analogous
  convergence of $Z_{n,\xi}$.

  Finally, by \eqref{eq:ZnxiZnxistarapprox}, the definition of $BL_1(\ell^\infty(\FF))$
   and Theorem \ref{theo:bootconv1}
  $$ \big| E_\xi g(Z_{n,\xi}^*)-    Eg(Z)\big| \le \big| E_\xi g(Z_{n,\xi}^*)-E_\xi
  g(Z_{n,\xi})\big| + \big|E_\xi  g(Z_{n,\xi})-    Eg(Z)\big|\;\longrightarrow\; 0
  $$
  in outer probability uniformly for all $g\in
  BL_1(\ell^\infty(\FF))$.
\end{proofof}

\begin{proofof} Theorem \ref{theo:extremo1}.\rm\quad
  The convergence of $(Z_n^{(h,\tilde h)}(f))_{f\in\FF}$ follows from Corollary 3.6(ii) and Remark 3.7(i) of Drees and Rootz\'{e}n (2010); see also Drees and Rootz\'{e}n (2015). To see this, check that
   Condition (D3) is fulfilled since for $\delta<1$
  \begin{eqnarray*}
    \lefteqn{\sup_{\varrho(f_D,f_{\tilde D})<\delta} \frac 1{r_n v_n} E\Big( \sum_{i=1}^{r_n}
   1_{D\Delta\tilde D}(X_{n,i}^{(h,\tilde h)})\Big)^2}\\
   & \le & \max_{\bar h\in\{0,\ldots,h_0\}}\sup_{\bar\varrho((A,B),(\tilde A,\tilde
   B))<\delta} \frac 1{r_n v_n} E\Big( \sum_{i=1}^{r_n}
   1_{(A\times B)\Delta (\tilde A\times \tilde
   B)}(X_i/a_k,X_{i+\bar h}/a_k)\Big)^2\\
   & \le & \sup_{0<t< \delta} u(t).
  \end{eqnarray*}
  Since, by ($\widetilde{\text{B3}}$), for $i>\max(h,\tilde h)$
  \begin{eqnarray*}
    P\big(X_{n,i+1}^{(h,\tilde h)}\ne 0\mid X_{n,1}^{(h,\tilde h)}\ne 0\big)
    & \le & \sum_{j,l\in\{0,h,\tilde h\}} P(X_{n,i+1+j}\ne 0\mid X_{n,1+l}\ne 0) \\
    & \le & \sum_{j,l\in\{0,h,\tilde h\}} s_n(i+j-l),
  \end{eqnarray*}
  $\big(X_{n,i}^{(h,\tilde h)}\big)_{1\le i\le n}$ satisfies the analog to ($\widetilde{\text{B3}}$) if $s_n(i)$ is replaced with
  $$ \tilde s_n(i) := \left\{ \begin{array}{l@{\quad}l}
     \sum_{j,l\in\{0,h,\tilde h\}} s_n(i+j-l), & i>\max(h,\tilde h),\\
     1, & i\le \max(h,\tilde h).
     \end{array}
     \right.
  $$

  Moreover,
  \eqref{eq:fourdimregvar} ensures that convergence (3.8) of Drees and Rootz\'{e}n (2010) holds, because
  \begin{eqnarray*}
   \frac 1{v_n^{(h,\tilde h)}} E\big(1_{A\times B\times\R^d}(X_{n,0}^{(h,\tilde h)}),1_{\tilde A\times\tilde
  B\times\R^d}(X_{n,i}^{(h,\tilde h)})\big) & = & \frac{k P\{a_k^{-1}(X_0,X_h,X_i,X_{i+
  h})\in A\times B\times \tilde A\times \tilde B\}}{k P\{a_k^{-1}(X_0,X_h,X_{\tilde h})\in \R^{3d}\setminus (-\infty,x_*)^{3d}\}}\\
  & \to & \frac{\nu_{(0,h,i,i+h)}(A\times B\times \tilde A\times \tilde
  B)}{\nu_{(0,h,\tilde h)}(\R^{3d}\setminus (-\infty,x_*)^{3d})}\\
  & =: & d_i \big(f_{A\times B\times\R^d}, f_{\tilde A\times\tilde
  B\times\R^d}\big),
  \end{eqnarray*}
  and likewise
  \begin{eqnarray*}
    \frac 1{v_n^{(h,\tilde h)}} E\big(1_{A\times B\times\R^d}(X_{n,0}^{(h,\tilde h)}),1_{\tilde A\times\R^d\times\tilde
  B}(X_{n,i}^{(h,\tilde h)})\big)
   & \to &
  \frac{\nu_{(0,h,i,i+\tilde h)}(A\times B\times \tilde A\times \tilde
  B)}{\nu_{(0,h,\tilde h)}(\R^{3d}\setminus (-\infty,x_*)^{3d})}\\
  & =: & d_i \big(f_{A\times B\times\R^d}, f_{\tilde A\times\R^d\times\tilde
  B}\big),\\
  \frac 1{v_n^{(h,\tilde h)}} E\big(1_{A\times \R^d\times B}(X_{n,0}^{(h,\tilde h)}),1_{\tilde A\times\tilde B\times\R^d}
  (X_{n,i}^{(h,\tilde h)})\big)
   & \to &
  \frac{\nu_{(0,\tilde h,i,i+h)}(A\times B\times \tilde A\times \tilde
  B)}{\nu_{(0,h,\tilde h)}(\R^{3d}\setminus (-\infty,x_*)^{3d})}\\
  & =: & d_i \big(f_{A\times \R^d\times B}, f_{\tilde A\times\tilde B\times\R^d}\big)\\
  \frac 1{v_n^{(h,\tilde h)}} E\big(1_{A\times \R^d\times B}(X_{n,0}^{(h,\tilde h)}),1_{\tilde A\times\R^d\times\tilde B}
  (X_{n,i}^{(h,\tilde h)})\big)
   & \to &
  \frac{\nu_{(0,\tilde h,i,i+\tilde h)}(A\times B\times \tilde A\times \tilde
  B)}{\nu_{(0,h,\tilde h)}(\R^{3d}\setminus (-\infty,x_*)^{3d})}\\
  & =: & d_i \big(f_{A\times \R^d\times B}, f_{\tilde A\times\R^d\times\tilde B}\big).
  \end{eqnarray*}
  Hence, by Drees and Rootz\'{e}n (2015), condition (C3) holds and $Z_n^{(h,\tilde h)}$ converges to a Gaussian process
   with the covariance function specified in formula (3.10) of Drees and Rootz\'{e}n (2010) in terms of the functions $d_i$.

  Since
  \begin{eqnarray*}
    \frac{v_n^{(h,\tilde h)}}{v_n} & = & \frac{k P\{a_k^{-1}(X_0,X_h,X_{\tilde h})\in \R^{3d}\setminus
    (-\infty,x_*)^{3d}\}}{k P\{a_k^{-1}X_0\in \R^{d}\setminus
    (-\infty,x_*)^{d}\}}\\
    & \to & \frac{\nu_{(0,h,\tilde h)}(\R^{3d}\setminus
    (-\infty,x_*)^{3d})}{\nu_0 (\R^d\setminus (-\infty,x_*)^d)},
  \end{eqnarray*}
  the convergence of $(\tilde Z_n(\bar h,A,B))_{\bar h\in\{h,\tilde
  h),(A,B)\in\CC}$ to a Gaussian process with covariance function $\tilde c$
  follows from the approximation (3.6) of Drees and
  Rootz\'{e}n (2010). Now the assertion is obvious.
\end{proofof}

\begin{proofof} condition (D6) in Example \ref{ex:orthants}.
\rm\quad
    For fixed $r\in\N$ define
  functions $f_D^{(r)}: \R^{2rd}\to\R$, $f_D^{(r)}(y_1,\ldots,y_r)
  :=\sum_{i=1}^r 1_D(y_i)$ with $D\in\{A\times B \mid
  (A,B)\in\CC\}=\{(x,\infty)\mid x=(x_1,\ldots,x_d)\in (x_*,\infty)^{2d}\}$. The
  subgraph of $f^{(r)}_{(x,\infty)}$ equals
  \begin{eqnarray*}
   \lefteqn{\big\{ (t,(y_1,\ldots,y_r))\in\R^{2rd+1} \mid
  t<f^{(r)}_{(x,\infty)}(y_1,\ldots, y_r)\big\}}\\
    & = & \bigcup_{j=0}^r
  (-\infty,j)\times \{y=(y_1,\ldots, y_r)\mid
  f^{(r)}_{(x,\infty)}(y)=j\} \\
    & =: & M_x.
  \end{eqnarray*}

  Consider some fixed set $S=\{(t^{(l)},(y_1^{(l)},\ldots,
  y_r^{(l)})) \mid 1\le l\le m\}$ of $m$ points in $\R^{2rd+1}$. If
  for  $x,\tilde x\in \R^{2d}$ the symmetric difference
  $(x,\infty)\Delta (\tilde x,\infty)$ does not contain any of the
  $y_i^{(l)}=(y_{i,1}^{(l)},\ldots,y_{i,2d}^{(l)})$, $1\le i\le r$, $1\le l\le m$, then the intersections
  $S\cap M_x$ and $S\cap M_{\tilde x}$ are identical. Since the
  hyperplanes $\{x\in\R^{2d}\mid x_j=y_{i,j}^{(l)}\}$, $1\le j\le
  2d$, $1\le i\le r$, $1\le l\le m$, divide $\R^{2d}$ into at most
  $(mr+1)^{2d}$ hypercubes and for $x,\tilde x$ belonging
   to the same hypercube $(x,\infty)\Delta (\tilde x,\infty)$ does not contain any of the $y_i^{(l)}$,
   the family $\CC$ can pick out at most
  $(mr+1)^{2d}$ different subsets of $S$. Hence it cannot shatter
  $S$ if $(mr+1)^{2d}<2^m$, which is fulfilled if $m\ge 3d\log r$
  and
  $r$ is sufficiently large.

  To sum up, so far we have shown that, for some $r_0\in\N$ and all $r\ge r_0$, the VC-index of $\FF^{(r)}
  :=\{f_{A\times B}^{(r)} \mid (A,B)\in\CC\}$ is less than $3d\log r$.
  By Theorem 2.6.7 of van der Vaart and Wellner (1996), we conclude
  that
  \begin{equation}  \label{eq:entropbd1}
    N\Big(\eps \big(\int (F^{(r)})^2\, dQ\big)^{1/2}, \FF^{(r)},
    L_2(Q)\Big) \le K_1 r^{K_2} \eps^{-K_3\log r}
  \end{equation}
  for all $\eps\in (0,1)$, all probability measures $Q$ on
  $\R^{2rd}$, and suitable universal constants $K_1$, $K_2$ and
  $K_3$ with $F^{(r)}(y):=\sum_{i=1}^r 1_{\R^{2d}\setminus(-\infty
  x_*]^{2d}}(y_i)$ denoting the envelope function of $\FF^{(r)}$.

  Next let $H(y) := \sum_{i=1}^r 1_{\{y_i\ne 0\}}$ for
  $y=(y_1,\ldots,y_r)\in\R^{2rd}$, $X_{n,i}^{(h)}:= (\tilde
  X_{n,i},\tilde X_{n,i+h})$ for $1\le i\le n$ and define
  independent copies $Y_{n,j}^{(h)*}$ of $Y_{n,j}^{(h)}:=
  (X_{n,i}^{(h)})_{(j-1)r_n<i\le jr_n}$, $1\le j\le m_n$. Consider
  the non-zero values of the $N_r := \sum_{j=1}^{m_n}1_{\{H(Y_{n,j}^{(h)*})\le r\}}$ of
  these blocks with at most $r$ non-zero $X_{n,i}^{(h)}$'s; if
  necessary, these are completed by zeros to obtain vectors $\bar
  Y_j := \big(Y_{n,j,i_1}^{(h)*}, \ldots, Y_{n,j,i_r}^{(h)*}\big)$,
  i.e. $Y_{n,j,i_l}^{(h)*}\ne 0$ for $1\le l\le H(Y_{n,j}^{(h)*})\le
  r$ and $Y_{n,j,i_l}^{(h)*}= 0$ for $H(Y_{n,j}^{(h)*})<l\le
  r$. Let
  $$Q_{n,r} := \frac 1{N_r} \sum_{j=1}^{m_n} \eps_{\bar Y_j}
  1_{\{H(Y_{n,j}^{(h)*})\le r\}},
  $$
  and consider the squared random $L_2$-distance
  \begin{eqnarray*}
    d_n^2(f_{(x,\infty)},f_{(\tilde x,\infty)})
      & = & \frac 1{n v_n^{(h,\tilde h)}} \sum_{j=1}^{m_n}
        \big(f_{(x,\infty)}(Y_{n,j}^{(h)*})  -  f_{(\tilde
        x,\infty)}(Y_{n,j}^{(h)*})\big)^2 \\
     & \le & \frac{N_r}{n v_n^{(h,\tilde h)}} \int (f_{(x,\infty)}-f_{(\tilde
        x,\infty)})^2 \, dQ_{n,r} + \frac1{n v_n^{(h,\tilde h)}}
         \sum_{j=1}^{m_n}     H^2(Y_{n,j}^{(h)*})1_{\{H(Y_{n,j}^{(h)*})>r\}}
    \end{eqnarray*}
    for all $r\in\N$. In particular,
    $$  d_n^2(f_{(x,\infty)},f_{(\tilde x,\infty)})
      \le  \frac{N_{R_{n,\eps}}}{n v_n^{(h,\tilde h)}} \int (f_{(x,\infty)}-f_{(\tilde
        x,\infty)})^2 \, dQ_{n,{R_{n,\eps}}} + \frac{\eps^2}2
    $$
  with
  $$ R_{n,\eps} := \max\bigg(\min \Big\{ r\in\N \;\Big|\; \frac1{n v_n^{(h,\tilde h)}}
         \sum_{j=1}^{m_n}
         H^2(Y_{n,j}^{(h)*})1_{\{H(Y_{n,j}^{(h)*})>r\}}<\frac{\eps^2}2\Big\},r_0\bigg),
  $$
  so that a ball with radius $\tilde\eps:=\big(nv_n^{(h,\tilde h)}/(2N_{R_{n,\eps}})\big)^{1/2} \eps$ w.r.t.\ $L_2(Q_{n,{R_{n,\eps}}})$ is contained in a ball with radius $\eps$ w.r.t.\ $d_n$.
  Note that
  $$ \int (F^{(r)})^2\, dQ_{R_{n,\eps}} \le \frac 1{N_{R_{n,\eps}}}
  \sum_{j=1}^{m_n} H^2(Y_{n,j}^{(h)*}) 1_{\{H(Y_{n,j}^{(h)*})\le
  R_{n,\eps}\}}.
  $$
  Hence, in view of \eqref{eq:entropbd1},  $\bar\FF$ (defined in \eqref{eq:defFbar}) can be covered by
  \begin{eqnarray*}
   \lefteqn{N\Big(\tilde\eps ,\FF^{(R_{n,\eps})},
  L_2(Q_{n,R_{n,\eps}})\Big)}\\
   & \le & K_1 R_{n,\eps}^{K_2} \bigg( \eps \Big(\frac{n v_n^{(h,\tilde
  h)}}{2\sum_{j=1}^{m_n} H^2(Y_{n,j}^{(h)*}) 1_{\{H(Y_{n,j}^{(h)*})\le
  R_{n,\eps}\}}}\Big)^{1/2}\bigg)^{-K_3 \log R_{n,\eps}}\\
   & \le & K_1 R_{n,\eps}^{K_2} \bigg( \frac{\eps}{R_{n,\eps}} \Big(\frac{2\sum_{j=1}^{m_n} 1_{\{Y_{n,j}^{(h)*}\ne 0\}}}  {n v_n^{(h,\tilde  h)}}\Big)^{-1/2}\bigg)^{-K_3 \log R_{n,\eps}}
  \end{eqnarray*}
  balls with radius $\eps$ w.r.t.\ $d_n$.

  Next observe that \eqref{eq:extremomomcond} implies $E\big(H^{2+\delta}(Y_{n,1}^{(h)*})\big)=O(r_n v_n)$:
  \begin{eqnarray}
    E\Big( \sum_{i=1}^{r_n} 1_{\{X_{n,i}^{(h)}\ne
    0\}}\Big)^{2+\delta}
    & \le & E\Big( 2\max_{l\in\{0,h\}}
        \sum_{i=1}^{r_n} 1_{\{\tilde X_{n,i+l}\ne
        0\}}\Big)^{2+\delta}  \nonumber\\
    & \le & 2^{2+\delta} \sum_{l\in\{0,h\}} E\Big(
        \sum_{i=1}^{r_n} 1_{\{\tilde X_{n,i+l}\ne
        0\}}\Big)^{2+\delta} \nonumber\\
    & = & 2^{3+\delta} E\Big(
        \sum_{i=1}^{r_n} 1_{\{X_i\not\in
        (-\infty,a_kx_*)^d\}}\Big)^{2+\delta}\nonumber \\
    & = & O(r_nv_n),  \label{eq:extremomomcond2}
  \end{eqnarray}
  where in the last but one line we have used the stationarity of
  the time series. Hence,
  $E\big(H(Y_{n,1}^{(h)*})\big)=r_nP\{a_k^{-1}(X_0,X_n)\not\in
  (-\infty,x_*)^{2d}\} =: r_n v_n^{(h)}=O(r_nv_n)$ implies that $r_n v_n^{(h)} =O\big(  P\{Y_{n,1}^{(h)*}\ne 0\}\big)$, because else
  $\limsup_{n\to\infty}E(H(Y_{n,1}^{(h)*})|Y_{n,1}^{(h)*}\ne
  0)=\infty$ and thus
  $$ \limsup_{n\to\infty}
  \frac{E(H^2(Y_{n,1}^{(h)*}))}{E(H(Y_{n,1}^{(h)*}))} =\limsup_{n\to\infty}
  \frac{E(H^2(Y_{n,1}^{(h)*})\mid Y_{n,1}^{(h)*}\ne
  0)}{E(H(Y_{n,1}^{(h)*})\mid Y_{n,1}^{(h)*}\ne
  0)} \ge \limsup_{n\to\infty} E(H(Y_{n,1}^{(h)*})\mid Y_{n,1}^{(h)*}\ne
  0) = \infty,
  $$
  in contradiction to \eqref{eq:extremomomcond2}. By Chebyshev's inequality,
  $$ P\Big\{ \sum_{j=1}^{m_n} 1_{\{Y_{n,j}^{(h)*}\ne 0\}}> 2 m_nP\{Y_{n,1}^{(h)*}\ne 0\}\Big\} \le \frac 1{m_n P\{Y_{n,1}^{(h)*}\ne 0\}} \to 0.
  $$
  Since $v_n$, $v_n^{(h)}$ and $v_n^{(h,\tilde h)}$ are all of the same order (by the regular variation of $(X_0,X_h,X_{\tilde h})$), we conclude that with
  probability tending to 1
  $$ N(\eps,\bar\FF,d_n) \le K_1 R_{n,\eps}^{K_2} (K_4 \eps/R_{n,\eps})^{-K_3
  \log R_{n,\eps}}.
  $$
  Finally, \eqref{eq:extremomomcond2} implies that to each $\eta>0$ there
  exist constants $M,\tau>0$ such that
  \begin{eqnarray*}
    \lefteqn{P\Big\{ \frac1{n v_n^{(h,\tilde h)}} \sum_{j=1}^{m_n}
         H^2(Y_{n,j}^{(h)*})
         1_{\{H(Y_{n,j}^{(h)*})>M\eps^{-(2+\tau)/\delta}\}}
         >\frac{\eps^2}2 \text{ for some } 0<\eps\le 1\Big\}}\\
   & \le & \sum_{l=0}^\infty P\Big\{ \frac1{n v_n^{(h,\tilde h)}} \sum_{j=1}^{m_n}
         H^2(Y_{n,j}^{(h)*})
         1_{\{H(Y_{n,j}^{(h)*})>M2^{l(2+\tau)/\delta}\}}
         >\frac{2^{-2(l+1)}}2\Big\} \\
   & \le & \sum_{l=0}^\infty 2^{2l+3} E\Big(\frac1{n v_n^{(h,\tilde h)}} \sum_{j=1}^{m_n} H^2(Y_{n,j}^{(h)*})
         1_{\{H(Y_{n,j}^{(h)*})>M2^{l(2+\tau)/\delta}\}}\Big) \\
   & \le & \sum_{l=0}^\infty \frac{2^{2l+3}}{n v_n^{(h,\tilde h)}} \cdot
        \frac{m_n E(H^{2+\delta}(Y_{n,1}^{(h)*}))}{(M
       2^{l(2+\tau)/\delta})^\delta} \\
   & \le & \frac{K_6}{M^\delta} \sum_{l=0}^\infty 2^{-2 l\tau}\\
   & < & \eta
  \end{eqnarray*}
  with $K_6$ denoting some universal constant.
  Hence $R_{n,\eps} \le M\eps^{-(2+\tau)/\delta}$ with probability
  greater than $1-\eta$, so that
  $$ \int_0^\xi \big(\log N(\eps,\bar\FF,d_n)\big)^{1/2}\, d\eps \le
  \int_0^\xi \big(K_7+K_8|\log\eps|+K_9(\log\eps)^2\big)^{1/2}\,
  d\eps
  $$
  tends to 0 as $\xi$ tends to 0, which proves condition (D6).
  Hence,
  under the additional assumptions of Theorem \ref{theo:extremo1},
  the process $\tilde Z_n$ converges.
\end{proofof}

\begin{proofof} Corollary \ref{corol:empextremo}.\rm\quad
 Check that
 \begin{eqnarray*}
  \hat\rho_{n,A,B}(h) & = & \frac{nP\{X_0\in a_kA,X_h\in a_kB\}+\sqrt{nv_n}\tilde Z_n(h,A,B)}{nP\{X_0\in a_kA\}
   +\sqrt{nv_n}\tilde Z_n(h,A,\R^d)}\\
   & = & \frac{\rho_{a_k,A,B}(h)+\sqrt{v_n/n}\tilde Z_n(h,A,B)/P\{X_0\in
   a_kA\}}{1+\sqrt{v_n/n}\tilde Z_n(h,A,\R^d)/P\{X_0\in   a_kA\}}\\
   & = & \rho_{a_k,A,B}(h)+\frac{\sqrt{v_n/n}}{P\{X_0\in
   a_kA\}}\cdot\frac{\tilde Z_n(h,A,B)-\rho_{a_k,A,B}(h)\tilde
   Z_n(h,A,\R^d)}{1+\sqrt{v_n/n}\tilde Z_n(h,A,\R^d)/P\{X_0\in
   a_kA\}}.
 \end{eqnarray*}
 Since by the regular variation of $X_0$
 $$ \frac{P\{X_0\in a_kA\}}{v_n}= \frac{P\{X_0\in
 a_kA\}}{P\{X_0\not\in (-\infty,x_*)^d\}} \to
 \frac{\nu_0(A)}{\nu_0\big(\R^d\setminus(-\infty,x_*)^d\big)},
 $$
 the first assertion is an immediate consequence of Theorem
 \ref{theo:extremo1} and the second follows from $v_n\sim
 \nu_0\big(\R^d\setminus(-\infty,x_*)^d\big)/k$.
\end{proofof}

\begin{proofof}  Corollary \ref{corol:empextremo2}.\rm\quad
  The first assertion follows from
  $ \hat{\hat\rho}_{n,A,B}(h)-\rho_{\hat a_k,A,B}(h) = $ \linebreak
  $\hat\rho_{n,(\hat a_k/a_k)A,(\hat a_k/a_k)B}(h)-\rho_{a_k,(\hat a_k/a_k)A,(\hat
  a_k/a_k)B}(h)$, the uniform convergence in
  \eqref{eq:empextremoconv} and the continuity of $\tilde
  Z(h,\cdot,\cdot)$ w.r.t.\ $\bar\varrho$.

  Under condition \eqref{eq:secordrho} we have by the extended
  regular variation of $\Phi_h$ and the consistency of $\hat a_k$
  \begin{eqnarray*}
    \rho_{\hat a_k,A,B}(h) -\rho_{a_k,A,B}(h) & = & (\Phi_h(\hat
    a_k)-\Phi_h(a_k)) \Psi_h(A,B)+ o(|\Phi_h(\hat a_k)|+|\Phi(a_k)|)
    \\
     & \le & |\Phi_h(a_k)| \Big|\frac{\Phi_h(\hat
     a_k)}{\Phi_h(a_k)}-1\Big| |\Psi_h(A,B)|+ o(|\Phi_h(a_k)|)\\
     & = & o_P((k/n)^{1/2}),
  \end{eqnarray*}
  which proves the second assertion.
\end{proofof}

\begin{proofof} Theorem \ref{theo:extremobootstrap}.\rm\quad
  For $ h\in\{0,\ldots, h_0\}$ and $ (A,B)\in\CC$, let
  \begin{eqnarray*}
   \tilde Z_{n,\xi}(h,A,B) &:= & \frac{\sqrt{nv_n^{(h,\tilde h)}}}{\sqrt{nv_n}}
   Z_{n,\xi}^{(h,\tilde h)}(f_{A\times B\times
   \R^d})\\
   & = &  \frac 1{\sqrt{nv_n}}  \sum_{j=1}^{m_n} (1+\xi_j)
  \sum_{i=1}^{r_n} \Big(1_{A\times B} \big( a_k^{-1}
  (X_{(j-1)r_n+i},X_{(j-1)r_n+i+h}\big)\\
  & &  \hspace*{3cm}{  }-P\big\{a_k^{-1}
  (X_{(j-1)r_n+i},X_{(j-1)r_n+i+h})\in A\times B\big\}\Big)
   \end{eqnarray*}
  with $Z_{n,\xi}^{(h,\tilde h)}$ denoting the multiplier process
  pertaining to $Z_n^{(h,\tilde h)}$ (cf.\ \eqref{eq:multiplierdef}).
  By Proposition \ref{prop:bootequi}, Theorem \ref {theo:bootconv1}
  and the proof of Theorem \ref{theo:extremo1} (in particular, the
  convergence of $v_n^{(h,\tilde h)}/v_n$)
  \begin{equation} \label{eq:extremoprocbootstrap}
   \sup_{g\in BL_1 (\ell^\infty(\{0,\ldots,h_0\}\times
  \CC) )} \big| E_\xi g(\tilde Z_{n,\xi}) -Eg(\tilde Z)\big|
  \;\longrightarrow\; 0
  \end{equation}
  in outer probability.

  Let
  $$ g_j(h,A,B) = f_{A\times B\times \R^d}(Y_{n,j}^{(h,\tilde h)} ) =\sum_{i=1}^{r_n} 1_{A\times B} \big( a_k^{-1}
  (X_{(j-1)r_n+i},X_{(j-1)r_n+i+h})\big).
  $$
  Recall from the proof of Theorem \ref{theo:extremo1} that
  $$ \hat\rho_{n,A,B}(h) =  \frac{\sum_{j=1}^{m_n}
    g_j(h,A,B)}{\sum_{j=1}^{m_n} g_j(h,A,\R^d)} +
    o_p((nv_n)^{-1/2})
  $$
  (cf.\ also Corollary 3.6 of Drees and Rootz\'{e}n, 2010).
  Thus
  \begin{eqnarray*}
    \lefteqn{R_{n,\xi}(h,A,B)}\\
     & = & \sqrt{nv_n} \bigg(\frac{\sum_{j=1}^{m_n} (1+\xi_j)
    g_j(h,A,B)}{\sum_{j=1}^{m_n} (1+\xi_j)    g_j(h,A,\R^d)} - \frac{\sum_{j=1}^{m_n}
    g_j(h,A,B)}{\sum_{j=1}^{m_n} g_j(h,A,\R^d)} +
    o_p((nv_n)^{-1/2})\bigg)\\
    & = & \sqrt{nv_n}\frac{ \sum_{j=1}^{m_n} \xi_j g_j(h,A,B)- \sum_{j=1}^{m_n} \xi_j g_j(h,A,\R^d)\cdot
    \sum_{j=1}^{m_n}g_j(h,A,B)/ \sum_{j=1}^{m_n} g_j(h,A,\R^d)}{\sum_{j=1}^{m_n} (1+\xi_j)    g_j(h,A,\R^d)}+
    o_p(1).
  \end{eqnarray*}
  Note that
  $$ \sum_{j=1}^{m_n} \xi_j g_j(h,A,B) = \sqrt{nv_n} \tilde
  Z_{n,\xi}(h,A,B)+ \sum_{j=1}^{m_n} \xi_j r_n P\{(X_0,X_h)/a_k\in
  A\times B\},
  $$
  where according to the central limit theorem and the regular variation of $(X_0,X_h)$
  the second term is of the order $O_P\big(m_n^{-1/2}
  r_nv_n\big)=O_P\big(\sqrt{nv_n}\sqrt{r_nv_n}\big)=o_P\big(\sqrt{nv_n}\big)$.
  Hence
  \begin{eqnarray*}
    R_{n,\xi}(h,A,B) & = & nv_n  \frac{\tilde
    Z_{n,\xi}(h,A,B)-\hat\rho_{n,A,B}(h)\tilde
    Z_{n,\xi}(h,A,\R^d)+o_P(1)}{m_nr_n P\{X_0/a_k\in A\} +
    O_P(\sqrt{nv_n})} +o_P(1)\\
     & = & \frac{\nu_0\big(\R^d\setminus
     (-\infty,x_*)^d\big)}{\nu_0(A)} \Big(\tilde
    Z_{n,\xi}(h,A,B)-\hat\rho_{n,A,B}(h)\tilde
    Z_{n,\xi}(h,A,\R^d)\Big) + o_P(1)
  \end{eqnarray*}
  uniformly for $h\in\{0,\ldots, h_0\}, (A,B)\in\CC$. In the
  last step we have used that  by the regular variation of $X_0$ and
  the definition of $v_n$
  $$ \frac{P\{X_0/a_k\in A\}}{v_n} \,\to\, \frac{\nu_0(A)}{\nu_0\big(\R^d\setminus
     (-\infty,x_*)^d\big)}
  $$
  where, by assumption, $\nu_0(A)$ is bounded away from 0. Now we
  can conclude \eqref{eq:extremobootstrap1} from
  \eqref{eq:extremoprocbootstrap} and \eqref{eq:empextremoconv}.

  Finally, notice that $\big| g(\hat R_{n,\xi})-g(R_{n,\xi})\big| \le
  \sup_{h\in\{0,\ldots,h_0\}, (A,B)\in\CC} \big|  R_{n,\xi}(h,(\hat a_k/a_k)A,(\hat
  a_k/a_k)B)-$ \linebreak
   $R_{n,\xi}(h,A,B)\big|\to 0$ in outer probability  for all $g\in BL_1\big(\ell^\infty(\{0,\ldots,h_0\}\times
  \CC)\big)$, because $\hat
  a_k/a_k\to 1$ and $R_{n,\xi}(h,\cdot,\cdot)$ is asymptotically
  equicontinuous w.r.t.\ $\bar\varrho$. Therefore,
  \eqref{eq:extremobootstrap2} is an immediate consequence of
  \eqref{eq:extremobootstrap1}.
\end{proofof}

\section*{Appendix}

The following conditions were used by Drees and Rootz\'{e}n (2010, 2015). For
the ease of reference, we use the same numbering as in these papers.
\\[1ex]
{\bf (B1)}\; \; \; \parbox[t]{14.9cm}{The rows $(X_{n,i})_{1\le i\le
n}$ are stationary, $ \ell_n = o(r_n),   \; \ell_n\to\infty,$ $r_n =
o(n)$, $r_nv_n \to 0,$  $nv_n \to \infty$}
\\[1ex]
{\bf (B2)} \; \; \;$\beta_{n,l_n} n/r_n \to 0. $ \\[1ex]
{\bf ($\widetilde{\text{B3}}$)} \; \; \;\parbox[t]{14.8cm}{For all $n\in\N$ and all $1\le i\le r_n$ there exists $s_n(i)\ge P(X_{n,i+1}\ne 0\mid X_{n,1}\ne 0)$ such that $s_\infty(i):=\lim_{n\to\infty} s_n(i)$ exists and $\lim_{n\to\infty} \sum_{i=1}^{r_n} s_n(i) =\sum_{i=1}^\infty s_\infty(i)<\infty$.}

Recall that for $y=(y_1,\ldots, y_r)$ and $l<r$ we define
$y^{(r-l)}:=(y_1,\ldots,y_{r-l})$. \\[1ex]
 {\bf (C1)}\;\;\; \parbox[t]{14.5cm}{For $\Delta_n(f)  :=
f(Y_n) - f(Y_n^{(r_n - \ell_n)})$
\begin{eqnarray}
      E\Big( (\Delta_n(f)- E\Delta_n(f))^2
      \Ind{|\Delta_n(f)- E\Delta_n(f)|\le
      \sqrt{nv_n}}\Big) & = & o(r_n v_n)   \label{necsuffbbsbcond1} \nonumber\\
     P \big\{|\Delta_n(f)- E\Delta_n(f)|>
     \sqrt{nv_n}\big\} & = & o(r_n/n) \label{necsuffbbsbcond2}
     \nonumber
   \end{eqnarray}
   for all $f\in\FF$.}
\\[1ex] {\bf(C2)}\; \; \; \parbox[t]{17cm}{
$E\Big( (f(Y_n) - E f(Y_n))^2\Ind{|f(Y_n) - E
    f(Y_n)|>\eps\sqrt{nv_n}}\Big) = o(r_n v_n),$ \quad $\forall\,
    \eps>0, f\in\FF.$}
\vspace*{2mm}
\\[1ex] {\bf (C3)}\; \; \;  $\displaystyle \frac 1{r_nv_n} Cov\big( f(Y_n),
g(Y_n)\big) \to c(f,g), \quad \forall\, f,g\in\FF.$\\[2ex]
{\bf (D1)}\;\;\; \parbox[t]{14.5cm}{The index set $\FF$ consists of
cluster functionals $f$ such that $E(f^2(Y_n))$ is finite for all
$n\geq 1$ and such that the envelope function
$$ F(x) := \sup_{f\in\FF} |f(x)| $$
is finite for all $x\in E_\cup$.} \vspace*{2mm}
\\[1ex]{\bf (D2)}  $$E^* \Big( F(Y_{n}) \Ind{F(Y_{n})>\eps\sqrt{nv_n}}\Big)  =
    o\big(r_n\sqrt{v_n/n}\big), \quad \forall\, \eps>0.$$
\vspace*{2mm}
\\[1ex]{\bf (D2')}\begin{equation*}
        E^* \Big( F^2(Y_{n}) \Ind{F(Y_{n})>\eps\sqrt{nv_n}}\Big)
        = o(r_nv_n), \quad\forall\, \eps>0.
    \end{equation*}
\\[1ex]{\bf (D3)}\;\;\;  \parbox[t]{14.9cm}{There exists a semi-metric $\rho$ on
$\FF$ such that $\FF$ is totally bounded (i.e., for all $\eps>0$ the
set $\FF$ can be covered by finitely many balls with radius $\eps$
w.r.t.\ $\rho$) such that  $$\lim_{\delta\downarrow 0}
\limsup_{n\to\infty}
    \sup_{f,g\in\FF, \; \rho(f,g)<\delta} \frac 1{r_nv_n}
    E(f(Y_{n})-g(Y_{n}))^2  = 0.$$}\\[1ex]

Finally, we consider different entropy conditions, which measure the
complexity of the family $\FF$. The {\em bracketing number}
$N_{[\cdot]}(\eps,\FF,L_2^n)$ is defined as the smallest number
$N_\eps$ such that for each $n\in\N$ there exists a partition
$(\FF_{n,k}^\eps)_{1\le k\le N_\eps}$ of $\FF$ such that
\begin{equation}
  E^* \sup_{f,g\in\FF_{n,k}^\eps} \big(f(Y_{n})-g(Y_{n})\big)^2
  \le \eps^2 r_n v_n, \quad \forall\, 1\le k\le N_\eps.
\end{equation}
 For a given semi-metric $d$ on $\FF$, the (metric) {\em covering number}
   $N(\eps,\FF,d)$ is the minimum number of balls with radius $\eps$
   w.r.t.\ $d$ needed to cover $\FF$. The condition (D6) bounds the rate of increase of
   $N(\eps,\FF,d_n)$ as $\eps$ tends to 0 for the random semi-metric
   $$ d_n(f,g) := \Big( \frac 1{nv_n} \sum_{j=1}^{m_n} \big(
   f(Y_{n,j}^*)-g(Y_{n,j}^*)\big)^2\Big)^{1/2},
   $$
   that is the $L_2$-semi-metric w.r.t.\ to empirical measure
   $(nv_n)^{-1}\sum_{j=1}^{m_n} \eps_{Y_{n,j}^*}$,
   where $Y_{n,j}^*$, $1\le j\le m_n$, are i.i.d.\ copies of $Y_{n,1}$.
    In    (D6') we instead use the supremum of all covering numbers $N(\eps,\FF,d_Q)$
    where
     $d_Q(f,g):= \big(\int (f-g)^2\, dQ\big)^{1/2}$ and $Q$
   ranges over the set of discrete probability measures
   $\mathcal{Q}$ on $(E_\cup,\mathbb{E}_\cup)$.\\[2ex]
 {\bf (D4)} \parbox[t]{15cm}{$$\lim_{\delta\downarrow 0}
\limsup_{n\to\infty} \int_0^\delta
     \sqrt{\log N_{[\cdot]}(\eps,\FF,L_2^n)}\, d\eps  =  0.$$}
\\[1ex]
   {\bf (D5)}\;\;\;  \parbox[t]{14.9cm}{For all $\delta>0$, $n\in\N$, $(e_i)_{1\le i\le
   \floor{m_n/2}} \in \{-1,0,1\}^{\floor{m_n/2}}$ and $k\in\{1,2\}$ the
   map $\sup_{f,g\in\FF, \rho(f,g)<\delta}$
   $\sum_{j=1}^{\floor{m_n/2}} e_j\big(
   f(Y_{n,j}^*)-g(Y_{n,j}^*)\big)^k$
      is measurable.}\\[1ex]
\\[1ex]
   {\bf (D6)} \begin{equation*}
    \lim_{\delta\downarrow 0} \limsup_{n\to\infty}P^*\Big\{ \int_0^\delta
    \sqrt{ \log N(\eps,\FF,d_n)}\, d\eps > \tau \Big\} = 0, \quad
    \forall \tau>0.
  \end{equation*}
\\[1ex]
   {\bf (D6')}\;\;\;  \parbox[t]{14.9cm}{The envelope function $F$ is
   measurable with $E( F^2(Y_{n}))=O(r_nv_n)$ and
 \begin{equation*}
    \int_0^1 \sup_{Q\in\mathcal{Q}} \sqrt{ \log
    N(\eps{\textstyle(\int F^2dQ)^{1/2}},\FF,d_Q)}\, d\eps < \infty.
  \end{equation*}}

\bigskip

\noindent {\bf Acknowledgement:} I thank Anja Jan{\ss}en for helpful
discussions about regular variation on general cones and for providing R-code to calculate the extremogram for $t$-GARCH time series. The financial
support by the German Research foundation DFG via the grant JA
2160/1 is gratefully acknowledged.

\bigskip

\noindent {\large\bf References}
   \smallskip

\parskip1.2ex plus0.2ex minus0.2ex

\rueck
  Basrak, B., and Segers, J. (2008). Regularly varying multivariate time series. Universit\'{e} Catholique de Louvain,
  Institut de Statistique discussion paper 0717.

\rueck
  Bingham, N.H., Goldie, C.M., and Teugels, J.L.\ (1987). {Regular
  Variation}. Cambridge University Press.



\rueck
  Das, B., Mitra, A., and  Resnick S.\ (2013). Living on the
  multi-dimensional edge: seeking hidden risks using regular
  variation.    {\em Adv.\ Appl.\ Probab.}
    {\bf 45}, 139--163.

\rueck
  Davis, R., and Mikosch, T. (2009). The extremogram: A correlogram for
extreme events. {\em Bernoulli} {\bf 15}, 977-–1009.

\rueck
  Davis, R., Mikosch, T., and Cribben, I. (2012). Towards estimating extremal serial dependence via the
bootstrapped extremogram. {\em J.\ Econometr.} {\bf 170},   142--152.



\rueck
 Drees, H. (2003). Extreme quantile estimation for dependent data with applications to
finance. {\em Bernoulli} {\bf 9}, 617--657.

\rueck
 Drees, H. (2011).  Bias correction for estimators of the extremal index. Preprint,    arXiv:1107.0935v1

\rueck
 Drees, H., and Rootz\'{e}n, H. (2010). Limit Theorems for Empirical Processes of Cluster Functionals.
  {\em Ann.\ Statist.} {\bf 38}, 2145--2186.

\rueck
 Drees, H., and Rootz\'{e}n, H. (2015). Correction note to ``Limit Theorems for Empirical Processes of Cluster Functionals''. Preprint, arXiv:1510.09090v1.

\rueck
 Drees, H., Segers, J., and Warcho\l, M. (2015). Statistics for Tail Processes of Markov Chains. {\em Extremes} {\bf 18}, 369--402.


\rueck
 Ehlert, A., Fiebig, U.-R., Jan{\ss}en, A., and Schlather, M. (2015). Joint extremal behavior of hidden and observable time series with applications to GARCH processes. {\em Extremes} {\bf 18}, 109--140.


\rueck
  Kosorok, M.R. (2003). Bootstraps of sums of independent but not identically distributed stochastic processes.
  {\em J.\ Multiv. Analysis} {\bf 84}, 299--318.

\rueck
  Meinguet, T., and Segers, J. (2010). Regularly varying time series in Banach spaces. Preprint, arXiv:1001.3262v1.



\rueck
  Petrov, V.V. (1995). {\em Limit Theorems of Probability Theory.} Oxford Science Publication.

\rueck
  Resnick, S.I. (2007). Heavy-tail Phenomena. Springer.


\rueck
  van der Vaart, A.W., and Wellner, J.A. (1996). {\em Weak Convergence
  and Empirical Processes.} Springer.

\end{document}